\documentclass[a4paper,12pt]{amsart}
\usepackage{amssymb,amsmath,a4wide,graphicx,stmaryrd,fullpage,setspace,microtype}
\usepackage[pagebackref=true]{hyperref}
\usepackage[all]{xy}
\title{The geometry of cubical and regular transition systems}

\author[P. Gaucher]{Philippe Gaucher}

\address{Laboratoire PPS  (CNRS UMR 7126)\\
  Univ Paris Diderot\\
  Sorbonne Paris Cit\'e\\ Case 7014\\ 75205 PARIS Cedex 13\\ France}

\urladdr{http://www.pps.univ-paris-diderot.fr/{\~{}}gaucher/} 



\swapnumbers

\newcommand{\C}{\mathcal{C}}
\newcommand{\D}{\mathcal{D}}
\newcommand{\K}{\mathcal{K}}
\newcommand{\LL}{\mathcal{A}}
\newcommand{\I}{\mathcal{I}}
\newcommand{\W}{\mathcal{W}}
\newcommand{\F}{\mathcal{F}}

\newcommand{\de}{\partial}
\newcommand{\p}{\times}

\newtheorem{thm}{Theorem}[section]
\newtheorem{prop}[thm]{Proposition}
\newtheorem{lem}[thm]{Lemma}
\newtheorem{cor}[thm]{Corollary}

\newtheorem{defn}[thm]{Definition}
\newtheorem{nota}[thm]{Notation}

\newcommand{\bd}{\begin{defn}}
\newcommand{\ed}{\end{defn}}
\newcommand{\bp}{\begin{prop}}
\newcommand{\ep}{\end{prop}}
\newcommand{\bth}{\begin{thm}}
\renewcommand{\eth}{\end{thm}}
\newcommand{\bpf}{\begin{proof}}
\newcommand{\epf}{\end{proof}}
\newcommand{\bc}{\begin{cor}}
\newcommand{\ec}{\end{cor}}

\newcommand{\fL}[1]{\ar@{->}[ll]_-{#1}}
\newcommand{\fR}[1]{\ar@{->}[rr]^-{#1}}
\newcommand{\fRr}[1]{\ar@{->}[rrr]^-{#1}}
\newcommand{\fD}[1]{\ar@{->}[dd]_-{#1}}
\newcommand{\fU}[1]{\ar@{->}[uu]^-{#1}}
\newcommand{\f}[2]{\ar@{->}[#1]|{#2}}
\newcommand{\ff}[2]{\ar@2{->}[#1]|{#2}}
\newcommand{\frr}[1]{\ar@{->}[rrrr]^-{#1}}

\newcommand{\fl}[1]{\ar@{->}[l]_-{#1}}
\newcommand{\fr}[1]{\ar@{->}[r]^-{#1}}
\newcommand{\fd}[1]{\ar@{->}[d]_-{#1}}
\newcommand{\fu}[1]{\ar@{->}[u]^-{#1}}

\newcommand{\iso}{\cong}

\renewcommand{\leq}{\leqslant}
\renewcommand{\geq}{\geqslant}
\newcommand{\dd}[1]{\uparrow\!\!{#1}\!\!\uparrow}

\def\cartesien{%
  \ar@{-}[]+R+<6pt,-2pt>;[]+RD+<6pt,-6pt>%
  \ar@{-}[]+D+<2pt,-6pt>;[]+RD+<6pt,-6pt>%
}
\def\cocartesien{%
  \ar@{-}[]+L+<-6pt,+2pt>;[]+LU+<-6pt,+6pt>%
  \ar@{-}[]+U+<-2pt,+6pt>;[]+LU+<-6pt,+6pt>%
}
\def\hocartesien{%
  \ar@{-}[]+R+<6pt,-2pt>;[]+RD+<6pt,-6pt>_{h}%
  \ar@{-}[]+D+<2pt,-6pt>;[]+RD+<6pt,-6pt>%
}
\def\hococartesien{%
  \ar@{-}[]+L+<-6pt,+2pt>;[]+LU+<-6pt,+6pt>_{h}%
  \ar@{-}[]+U+<-2pt,+6pt>;[]+LU+<-6pt,+6pt>%
}

\newcommand{\brm}[1]{\rm{\mathbf{#1}}}

\newcommand{\set}{{\brm{Set}}}

\DeclareMathOperator{\id}{Id}

\DeclareMathOperator{\Mor}{Mor}

\newcommand{\liminj}{\varinjlim}

\newcommand{\wts}{\mathcal{W\!T\!S}}
\newcommand{\cts}{\mathcal{C\!T\!S}}
\newcommand{\rts}{\mathcal{R\!T\!S}}

\DeclareMathOperator{\cub}{\underline{Cub}}

\DeclareMathOperator{\dom}{dom}
\DeclareMathOperator{\bl}{\brm{\underline{L}}}

\makeatletter
\def\varholim@#1#2{%
  \vtop{\m@th\ialign{##\cr
    \hfil$#1\operator@font holim$\hfil\cr
    \noalign{\nointerlineskip\kern1.5\ex@}#2\cr
    \noalign{\nointerlineskip\kern-\ex@}\cr}}%
}
\def\holimproj{%
  \mathop{\mathpalette\varholim@{\leftarrowfill@\textstyle}}\nmlimits@
}
\def\holiminj{%
  \mathop{\mathpalette\varholim@{\rightarrowfill@\textstyle}}\nmlimits@
}
\makeatother

\DeclareMathOperator{\cell}{{\brm{cell}}}

\DeclareMathOperator{\inj}{{\brm{inj}}}
\newcommand{\ddownarrow}{{\downarrow}}

\DeclareMathOperator{\cyl}{{Cyl}}
\DeclareMathOperator{\cocyl}{{Path}}
\DeclareMathOperator{\CSA}{CSA}

\begin{document}

\begin{abstract}
  Il existe des syst\`emes de transitions cubiques contenant des cubes
  ayant un nombre arbitrairement grand de faces. Un syst\`eme de
  transition r\'egulier est un syst\`eme de transitions cubique tel
  que tout cube a le bon nombre de faces. Les propri\'et\'es
  cat\'egoriques et homotopiques des syst\`emes de transitions
  r\'eguliers sont similaires \`a celles des cubiques. On donne une
  description combinatoire compl\`ete des objets fibrants dans les cas
  cubiques et r\'eguliers. Un des deux appendices contient un lemme
  ind\'ependant sur la restriction d'une adjonction \`a une
  sous-cat\'egorie r\'eflective pleine.

  There exist cubical transition systems containing cubes having an
  arbitrarily large number of faces. A regular transition system is a
  cubical transition system such that each cube has the good number of
  faces. The categorical and homotopical results of regular transition
  systems are very similar to the ones of cubical ones.  A complete
  combinatorial description of fibrant cubical and regular transition
  systems is given. One of the two appendices contains a general lemma
  of independant interest about the restriction of an adjunction to a
  full reflective subcategory.
\end{abstract}

\maketitle

\tableofcontents

\section{Introduction}

\subsection*{Presentation}
\begin{figure}
\[
\xymatrix{
& \beta \ar@{->}[rd]^{\{b\}}&\\
\alpha\ar@{->}[ru]^{\{a\}}\ar@{->}[rd]_{\{b\}} & \{a,b\} & \delta\\
&\gamma\ar@{->}[ru]_{\{a\}}&}
\]
\caption{$a|| b$ : Concurrent execution of $a$ and $b$}
\label{concab}
\end{figure}
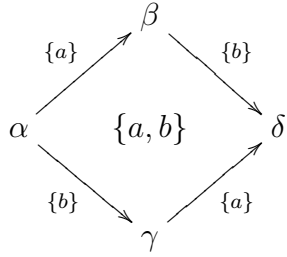

The purpose of Cattani-Sassone's notion of higher dimensional
  transition system introduced in \cite{MR1461821} is to model the
  concurrent execution of $n$ actions by a transition between two
  states labelled by a multiset $\{u_1,\dots,u_n\}$ of actions. A
  multiset is a set with a possible repetition of its elements:
  e.g. $\{u\}$ is not equal to $\{u,u\}$. A higher dimensional
  transition system for Cattani and Sassone consists of a set of
  states $S$, a set of actions $L$, a set of labels $\Sigma$ together
  with a labelling map $\mu:L\to \Sigma$, and a set of tuples
  $(\alpha,T,\beta)$ of transitions where $\alpha$ and $\beta$ are two
  states and $T$ is a multiset of actions. All these data have to
  satisfy several axioms which are detailed in the original paper
  \cite{MR1461821}. The higher dimensional transition $a||b$ depicted
  by Figure~\ref{concab} consists of the transitions
  $(\alpha,\{a\},\beta)$, $(\beta,\{b\},\delta)$,
  $(\alpha,\{b\},\gamma)$, $(\gamma,\{a\},\delta)$ and
  $(\alpha,\{a,b\},\delta)$. The labelling map is the identity map.
  Note that with $a=b$, we would get the $2$-dimensional transition
  $(\alpha,\{a,a\},\delta)$ which is not equal to the $1$-dimensional
  transition $(\alpha,\{a\},\delta)$. The latter actually does not
  exist in Figure~\ref{concab}. Indeed, the only $1$-dimensional
  transitions labelled by the multiset $\{a\}$ are
  $(\alpha,\{a\},\beta)$ and $(\gamma,\{a\},\delta)$.

In \cite{hdts}, Cattani-Sassone's notion is reworded in a more
convenient mathematical setting by introducing the notion of weak
transition system.  In this new setting, the transition
$(\alpha,\{a,b\},\delta)$ is represented by the tuple
$(\alpha,a,b,\delta)$. The set of transitions has therefore to satisfy
the Multiset axiom (here: if the tuple $(\alpha,a,b,\delta)$ is a
transition, then the tuple $(\alpha,b,a,\delta)$ has to be a
transition as well) and the Composition axiom which is a topological
version (in the sense of topological functors) of Cattani-Sassone's
interleaving axioms. The Composition axiom is called the Coherence
axiom in \cite{hdts}. The terminology is changed in the next paper
\cite{cubicalhdts} because this axiom behaves a little bit like a
partial $5$-ary composition in the proofs~\footnote{In the nLab
    page devoted to higher dimensional transition systems, T. Porter
    uses the terminology ``patching axiom'', which is quite a good
    idea too.}.  For example, the Composition axiom is the key axiom
for interpreting the higher dimensional transition system modeling the
$n$-cube as the free object generated by a ``pure'' $n$-dimensional
transition (this weak transition system consists of two states and a
$n$-dimensional transition going from one state to the other one)
\cite[Theorem~5.6]{hdts}. Indeed, the free compositions generated by
the Composition axiom generate all transitions of lower dimension
between the intermediate states (i.e. with a source different from the
initial state and a target different from the final state) .  Weak
transition systems assemble into a locally finitely presentable
category $\wts$ such that the forgetful functor forgetting the
transitions, and keeping the states and the actions, is topological in
the sense of \cite[Definition~21.1]{topologicalcat}.

The full coreflective subcategory $\cts$ of cubical transition systems
was then introduced in \cite{cubicalhdts}. They consist of the weak
transition systems which are equal to the union of their subcubes. It
was preferred to the full coreflective category of $\wts$ of colimits
of cubes because the latter does not contain the boundary of a
$2$-cube.  The category $\cts$ is sufficient to describe the path
spaces of all process algebras for any synchronization algebra because
their path spaces are colimits of cubes and because all colimits of
cubes are unions of cubes. Indeed, the weak transition system
associated with a process algebra is obtained by starting from a
labelled precubical set using the method described in
\cite{ccsprecub}, and by taking the free symmetric labelled precubical
set generated by it \cite{symcub}, and then by applying the
colimit-preserving realization functor from labelled symmetric
precubical sets to weak transition systems constructed in \cite{hdts}.

However, the notion of cubical transition system is still too
general. Indeed, a $n$-dimensional transition in a cubical transition
system may have an arbitrarily large number of faces in each
dimension.  Here is a simple example of a $2$-transition $X$ with $2n+2$
edges for an arbitrarily large integer $n\geq 1$:
\begin{itemize}
\item the set of states is
$\{I,F,a,b_1,\dots,b_n\}$
\item the set of actions is $\{u,v\}$ with
$\mu(u)\neq \mu(v)$ ($\mu$ denotes the labelling map)
\item the transitions are the tuples $(I,u,v,F)$,
$(I,v,u,F)$, $(I,u,a)$, $(a,v,F)$, $(I,v,b_i)$ and $(b_i,u,F)$ for
$i\leq 1\leq n$.
\end{itemize}
The weak transition system above is cubical because it is the union,
for $1\leq i \leq n$, of the $2$-cubes $Z_i$ having the set of
vertices $\{I,F,a,b_i\}$, the set of actions $\{u,v\}$ and the set of
six transitions
$\{(I,u,v,F),(I,v,u,F),(I,u,a),(a,v,F),(I,v,b_i),(b_i,u,F)\}$. To
avoid such a behavior, it suffices to replace the Intermediate state
axiom by the Unique intermediate state axiom, also called CSA2 (see
Definition~\ref{def_HDTS}). The latter axiom is already introduced in
\cite{hdts} to formalize Cattani-Sassone's notion of higher
dimensional transition systems in the setting of weak transition
systems.  We obtain a full reflective subcategory $\rts$ of that of
cubical transition systems whose objects are called the \emph{regular
  transition systems}.  Roughly speaking, a regular transition system
is a Cattani-Sassone transition system which does not necessarily
satisfy CSA1 (see Definition~\ref{csa1}). There is the chain of
functors
\[\rts \subset_{\hbox{\tiny reflective}} \cts \subset_{\hbox{\tiny coreflective}} \wts \stackrel{\omega}\longrightarrow_{\hbox{\tiny topological}} \set^{\{s\}\cup \Sigma}\]
where $\omega$ is the topological functor towards a power of the
category of sets forgetting the transitions: $s$ denotes the sort of
states and each element $x$ of the set of labels $\Sigma$ denotes the
sort of actions labelled by $x$. With the notations above, one
has \[\omega(a|| b) = (\{\alpha,\beta,\gamma,\delta\},\{a\},
\{b\})\] since the labelling map is the identity map. One has
\[\omega(X)=(\{I,F,a,b_1,\dots,b_n\}, \{u\}, \{v\}) \] since 
$\mu(u)\neq \mu(v)$.

Note that none of the categories of colimits of cubes and of regular
transition systems is included in the other one: see the final comment
of Section~\ref{reg-def}.

This paper is devoted to the geometric properties of regular
transition systems and to their relationship with cubical ones. Their
study requires the use of the whole chain of functors above which is
the composite of a right adjoint followed by a left adjoint followed
by a topological functor.  Despite the fact that colimits are
different in $\rts$ and in $\cts$, the main results are very similar
to the ones obtained for cubical transition systems in
\cite{cubicalhdts}. We will therefore follow the plan of
\cite{cubicalhdts}. The left determined model structure with respect
to the cofibrations of cubical transition systems between regular ones
is proved to exist. It is proved that the Bousfield localization by
the cubification functor is the model structure having the same class
of cofibrations and the fibrant objects are the regular transitions
systems such that for any transition $(\alpha,u_1,\dots,u_n,\beta)$,
the tuple $(\alpha,v_1,\dots,v_n,\beta)$ is a transition if
$\mu(u_i)=\mu(v_i)$ for $1\leq i \leq n$.  The homotopical structure
of this Bousfield localization will be completely elucidated. Roughly
speaking, after identifying each action of a regular transition system
with its label and after removing all non-discernable higher
dimensional transitions, two regular transition systems are weakly
equivalent if and only if they are isomorphic.

\subsection*{Outline of the paper}

Section~\ref{reg-def} introduces all definitions of higher dimensional
transition systems used in this paper: weak, cubical, regular. It
starts with the notion of \emph{regular transition system}
(Definition~\ref{def_HDTS}), and then by removing some axioms, the
notions of cubical transition system and of weak transition system are
recalled. This section does not contain anything new, except the
notion of regular transition system.
Section~\ref{section-cubical-lift} is a technical section which
provides a sufficient condition for an $\omega$-final lift of cubical
transition systems to be cubical (Theorem~\ref{cubical-lift}). This
result is used in the construction of several cubical transition
systems.  Section~\ref{cat-rts} deals with the most elementary
properties of regular transition systems.  The reflection
$\CSA_2:\cts\to \rts$ is studied.  The definition of the cubification
functor is recalled and its relationship with regular transition
systems is explained.  Section~\ref{homotopy-rts} establishes the
existence of the left determined model structure of regular transition
systems. The weak equivalences of this model structure are completely
characterized.  The Bousfield localization of the left determined
model category of regular transition systems by the cubification
functor is studied and completely elucidated in
Section~\ref{Bousfield-cub}. The comparison with cubical transition
systems is discussed there.  The proof of Theorem~\ref{compare-cellS}
is postponed to Section~\ref{cellS} to not overload
Section~\ref{Bousfield-cub}.  Finally, Section~\ref{fibrant-char}
completely characterizes the fibrant cubical and regular transition
systems in the Bousfield localizations by the cubification functor.
Section~\ref{restriction-adjunction-reflective-subcat} is a
categorical lemma of independant interest providing a easy way to
restrict an adjunction to a full reflective subcategory.

\subsection*{Prerequisites and notations}

All categories are locally small. The set of maps in a category $\K$
from $X$ to $Y$ is denoted by $\K(X,Y)$.  The initial (final resp.)
object, if it exists, is always denoted by $\varnothing$ ($\mathbf{1}$
resp.). The identity of an object $X$ is denoted by $\id_X$.  A
subcategory is always
isomorphism-closed.  
We refer to \cite{MR95j:18001} for locally presentable categories, to
\cite{MR2506258} for combinatorial model categories, and to
\cite{topologicalcat} for topological categories, i.e. categories
equipped with a topological functor towards a power of the category of
sets.  We refer to \cite{MR99h:55031} and to \cite{ref_model2} for
model categories. For general facts about weak factorization systems,
see also \cite{ideeloc}. The reading of the first part of
\cite{MOPHD}, published in \cite{MO}, is recommended for any reference
about good, cartesian, and very good cylinders.

\section{Regular higher dimensional transition systems}
\label{reg-def}

This section does not contain anything new, except the notion of
\emph{regular transition system}. It collects definitions and facts
about the various notions of transition systems which were expounded
in the previous papers of this series \cite{hdts} and
\cite{cubicalhdts}. To keep this section concise, the definition of a
regular transition system is given first, and then by removing some
axioms, the definitions of a cubical transition system and of a weak
transition system are recalled. It is necessary to recall all these
definitions because most of the proofs of this paper make use of the
whole chain of functors
\[\rts \subset_{\hbox{\tiny reflective}} \cts \subset_{\hbox{\tiny coreflective}} \wts \stackrel{\omega}\longrightarrow_{\hbox{\tiny topological}} \set^{\{s\}\cup \Sigma}\]
where $\set$ is the category of sets.

\begin{nota} A nonempty set of {\rm labels} $\Sigma$ is fixed.  \end{nota}

\bd\label{def_HDTS} A {\rm regular higher dimensional transition
  system} consists of a triple \[X=(S,\mu:L\rightarrow
\Sigma,T=\bigcup_{n\geq 1}T_n)\] where $S$ is a set of {\rm states},
where $L$ is a set of {\rm actions}, where $\mu:L\rightarrow \Sigma$
is a set map called the {\rm labelling map}, and finally where
$T_n\subset S\p L^n\p S$ for $n \geq 1$ is a set of {\rm
  $n$-transitions} or {\rm $n$-dimensional transitions} such that one
has:
\begin{itemize}
\item (All actions are used) For every $u\in L$, there is a
  $1$-transition $(\alpha,u,\beta)$.
\item (Multiset axiom) For every permutation $\sigma$ of
  $\{1,\dots,n\}$ with $n\geq 2$, if the tuple
  $(\alpha,u_1,\dots,u_n,\beta)$ is a transition, then the tuple
  $(\alpha,u_{\sigma(1)}, \dots, u_{\sigma(n)}, \beta)$ is a
  transition as well.
\item (Composition axiom~\footnote{This axiom is called the Coherence
    axiom in \cite{hdts} and \cite{cubicalhdts}.}) For every
  $(n+2)$-tuple $(\alpha,u_1,\dots,u_n,\beta)$ with $n\geq 3$, for
  every $p,q\geq 1$ with $p+q<n$, if the five tuples $(\alpha,u_1,
  \dots, u_n, \beta)$, $(\alpha,u_1, \dots, u_p, \nu_1)$, $(\nu_1,
  u_{p+1}, \dots, u_n, \beta)$, $(\alpha, u_1, \dots, u_{p+q}, \nu_2)$
  and $(\nu_2, u_{p+q+1}, \dots, u_n, \beta)$ are transiti\-ons, then
  the $(q+2)$-tuple $(\nu_1, u_{p+1}, \dots, u_{p+q}, \nu_2)$ is a
  transition as well.
\item (Unique intermediate state axiom or CSA2)~\footnote{This axiom
    is also called CSA2 in \cite{hdts}}.  For every
  $n\geq 2$, every $p$ with $1\leq p<n$ and every transition
  $(\alpha,u_1,\dots,u_n,\beta)$ of $X$, there exists a unique state
  $\nu$ such that both $(\alpha,u_1,\dots,u_p,\nu)$ and
  $(\nu,u_{p+1},\dots,u_n,\beta)$ are transitions.
\end{itemize}
A map of regular transition systems
\[f:(S,\mu : L \rightarrow \Sigma,(T_n)_{n\geq 1}) \rightarrow
(S',\mu' : L' \rightarrow \Sigma ,(T'_n)_{n\geq 1})\] consists of a
set map $f_0: S \rightarrow S'$, a commutative square
\[
\xymatrix{
  L \ar@{->}[r]^-{\mu} \ar@{->}[d]_-{\widetilde{f}}& \Sigma \ar@{=}[d]\\
  L' \ar@{->}[r]_-{\mu'} & \Sigma}
\] 
such that if $(\alpha,u_1,\dots,u_n,\beta)$ is a transition, then
$(f_0(\alpha),\widetilde{f}(u_1),\dots,\widetilde{f}(u_n),f_0(\beta))$
is a transition. The corresponding category is denoted by $\rts$.  The
$n$-transition $(\alpha,u_1,\dots,u_n,\beta)$ is also called a {\rm
  transition from $\alpha$ to $\beta$}.  The maps $f_0$ and
$\widetilde{f}$ will be also denoted by $f$. \ed

\begin{nota} The labelling map from the set of actions to the set of
  labels will be very often denoted by $\mu$. The set of states of a
  regular transition system $X$ is denoted by $X^0$. \end{nota}

The category $\rts$ of regular higher dimensional transition systems
is a full subcategory of the category of \emph{cubical transition
  systems} $\cts$ introduced in \cite{cubicalhdts}. By definition, a
cubical transition system satisfies all axioms of higher dimensional
transition system but one: the Unique intermediate state axiom is
replaced by the Intermediate state axiom, the state $\nu$ is not
necessarily unique anymore.  The category $\cts$ is a full subcategory
of the category of \emph{weak transition systems} $\wts$ introduced in
\cite{hdts}. By definition, a weak transition system satisfies all
axioms of regular transition systems but two: the Unique intermediate
state axiom is removed and an action is not necessarily used. Weak
transition system is the ``minimal'' definition: the multiset axiom is
indeed required to ensure that the concurrent execution of $n$ actions
does not depend on the order of the labelling, and the composition
axiom is required (even if its use is often hidden) e.g. to ensure
that labelled $n$-cubes are free objects (e.g. see the proof of
\cite[Theorem~5.6]{hdts}).  One has the inclusions of full
subcategories $\rts\subset \cts \subset \wts$. The inclusion
$\rts\subset \cts$ is strict since the introduction gives an example
of cubical transition system which is not regular.  The situation is
summarized in Table~\ref{alldefs}. Let us recall now the definition of
CSA1 for this sequence of definitions to be complete:

\bd \label{csa1} \cite[Definition~4.1~(2)]{hdts} and
\cite[Definition~7.1]{cubicalhdts} A cubical transition system
satisfies the {\rm First Cattani-Sassone axiom} (CSA1) if for every
transition $(\alpha,u,\beta)$ and $(\alpha,u',\beta)$ such that the
actions $u$ and $u'$ have the same label in $\Sigma$, one has $u =
u'$. \ed

\begin{table}[h]
\begin{tabular}{|l|c|c|c|c|}
\hline
& Cattani-Sassone & Regular& Cubical& Weak  \\
\hline
Multiset axiom & yes &  yes & yes & yes \\
\hline
Composition axiom & yes & yes & yes & yes \\
\hline
All actions used & yes &  yes & yes & no \\
\hline
Intermediate state axiom & yes & yes & yes & no\\
\hline
Unique Intermediate state axiom  & yes & yes & no & no \\
\hline
CSA1 & yes & no & no & no \\
\hline
\end{tabular}
\caption{Summary of all variants of transition systems.}
\label{alldefs}
\end{table}

The category $\wts$ is locally finitely presentable and the
functor \[\omega : \wts \longrightarrow \set^{\{s\}\cup \Sigma}\]
taking the weak higher dimensional transition system $(S,\mu : L
\rightarrow \Sigma,(T_n)_{n\geq 1})$ to the $(\{s\}\cup \Sigma)$-tuple
of sets $(S,(\mu^{-1}(x))_{x\in \Sigma}) \in \set^{\{s\}\cup \Sigma}$
is topological by \cite[Theorem~3.4]{hdts}. 

Let us recall that the paradigm of \emph{topological functor} is the
underlying set functor from the category of general topological spaces
to that of sets. The notion of topological functor is a generalization
of the notions of initial and final topologies
\cite{topologicalcat}. More precisely, a functor $\omega: \C
\rightarrow \D$ is \emph{topological} if each cone $(f_i:X \rightarrow
\omega A_i)_{i\in I}$ where $I$ is a class has a unique
$\omega$-initial lift (the \emph{initial structure})
$(\overline{f}_i:A \rightarrow A_i)_{i\in I}$, i.e.: 1) $\omega A= X$
and $\omega \overline{f}_i = f_i$ for each $i\in I$; 2) given $h:
\omega B \rightarrow X$ with $f_i h= \omega \overline{h}_i$,
$\overline{h}_i:B\rightarrow A_i$ for each $i\in I$, then $h=\omega
\overline{h}$ for a unique $\overline{h} : B \rightarrow A$.
Topological functors can be characterized as functors such that each
cocone $(f_i:\omega A_i \rightarrow X)_{i\in I}$ where $I$ is a class
has a unique $\omega$-final lift (the \emph{final structure})
$\overline{f}_i:A_i \rightarrow A$, i.e.: 1) $\omega A= X$ and $\omega
\overline{f}_i = f_i$ for each $i\in I$; 2) given $h: X \rightarrow
\omega B$ with $hf_i = \omega \overline{h}_i$, $\overline{h}_i:A_i
\rightarrow B$ for each $i\in I$, then $h=\omega \overline{h}$ for a
unique $\overline{h} : A \rightarrow B$.  A limit (resp. colimit) in
$\C$ is calculated by taking the limit (resp. colimit) in $\D$, and by
endowing it with the initial (resp. final) structure.  In particular,
a topological functor is faithful and it creates all limits and
colimits.

The category $\cts$ is a full coreflective locally finitely
presentable subcategory of $\wts$ by
\cite[Corollary~3.15]{cubicalhdts}. The composite functor
\[\cts \subset \wts \stackrel{\omega}\longrightarrow \set^{\{s\}\cup
  \Sigma}\] is faithful and colimit-preserving.

The inclusion $\cts\subset \wts$ is strict. Here are two families of
examples of weak transition systems which are not cubical:
\begin{enumerate}
\item The weak transition system $\underline{x} = (\varnothing, \{x\}
  \subset \Sigma, \varnothing)$ for $x \in \Sigma$ is not cubical
  because the action $x$ is not used. 
\item Let $n \geq 0$. Let $x_1,\dots,x_n \in \Sigma$. The \emph{pure
    $n$-transition} $C_n[x_1,\dots,x_n]^{ext}$ is the weak transition
  system with the set of states $\{0_n,1_n\}$, with the set of actions
  \[\{(x_1,1), \dots, (x_n,n)\}\] and with the transitions all
  $(n+2)$-tuples $(0_n,(x_{\sigma(1)},\sigma(1)), \dots,
  (x_{\sigma(n)},\sigma(n)),1_n)$ for $\sigma$ running over the set of
  permutations of the set $\{1,\dots ,n\}$. It is not cubical for $n>
  1$ because it does not contain any $1$-transition. Intuitively, the
  pure transition is a cube without faces of lower dimension.
\end{enumerate}

We give now some important examples of regular transition systems. In
each of the following examples, the axioms of regular higher
dimensional transition systems are satisfied for trivial reasons.

\begin{nota} For $n\geq 1$, let $0_n = (0,\dots,0)$ ($n$-times) and
  $1_n = (1,\dots,1)$ ($n$-times). By convention, let
  $0_0=1_0=()$. \end{nota}

\begin{enumerate}
\item Every set $X$ may be identified with the cubical transition
  system having the set of states $X$, with no actions and no
  transitions.
\item For every $x\in \Sigma$, let us denote by $\dd{x}$ the cubical
  transition system with four states $\{1,2,3,4\}$, one action $x$ and
  two transitions $(1,x,2)$ and $(3,x,4)$. The cubical transition
  system $\dd{x}$ is called the \emph{double transition (labelled by
    $x$)} where $x\in \Sigma$.
\end{enumerate}

Let us introduce now the cubical transition system corresponding to
the labelled $n$-cube.

\bp \label{cas_cube} \cite[Proposition~5.2]{hdts} Let $n\geq 0$ and
$x_1,\dots,x_n\in \Sigma$. Let $T_d\subset \{0,1\}^n \p
\{(x_1,1),\dots,(x_n,n)\}^d \p \{0,1\}^n$ (with $d\geq 1$) be the
subset of $(d+2)$-tuples
\[((\epsilon_1,\dots,\epsilon_n), (x_{i_1},i_1),\dots,(x_{i_d},i_d),
(\epsilon'_1,\dots,\epsilon'_n))\] such that
\begin{itemize}
\item $i_m = i_n$ implies $m = n$, i.e. there are no repetitions in the
  list $(x_{i_1},i_1),\dots,(x_{i_d},i_d)$
\item for all $i$, $\epsilon_i\leq \epsilon'_i$
\item $\epsilon_i\neq \epsilon'_i$ if and only if
  $i\in\{i_1,\dots,i_d\}$. 
\end{itemize}
Let $\mu : \{(x_1,1),\dots,(x_n,n)\} \rightarrow \Sigma$ be the set
map defined by $\mu(x_i,i) = x_i$. Then \[C_n[x_1,\dots,x_n] =
(\{0,1\}^n,\mu : \{(x_1,1),\dots,(x_n,n)\}\rightarrow
\Sigma,(T_d)_{d\geq 1})\] is a well-defined cubical transition system
called the {\rm $n$-cube}. \ep

The $n$-cubes $C_n[x_1,\dots,x_n]$ for all $n\geq 0$ and all
$x_1,\dots,x_n\in \Sigma$ are regular by \cite[Proposition~5.2]{hdts}
and \cite[Proposition~4.6]{hdts}.  For $n = 0$, $C_0[]$, also denoted
by $C_0$, is nothing else but the one-state higher dimensional
transition system $(\{()\},\mu:\varnothing \rightarrow
\Sigma,\varnothing)$.  

The category $\cts$ is a small-injectivity class of $\wts$ by
\cite[Theorem~3.6]{cubicalhdts}: more precisely being cubical is
equivalent to being injective with respect to the set of inclusions
$C_n[x_1,\dots,x_n]^{ext} \subset C_n[x_1,\dots,x_n]$ and
$\underline{x_1} \subset C_1[x_1]$ for all $n\geq 0$ and all
$x_1,\dots, x_n \in \Sigma$. Note that the composition axiom plays a
central role in this result.

Finally, let us notice that there is the isomorphism of weak
transition systems
\[\dd{x} \iso \liminj \left(C_1[x] \leftarrow \underline{x} \rightarrow C_1[x]\right)\] 
for any label $x$ of $\Sigma$, the colimit being taken in $\wts$. The
double transition $\dd{x}$ is an example of cubical transition system,
and even of regular transition system, which is not a colimit of
cubes. Another example of regular transition system which is not a
colimit of cubes is the boundary of a labelled $2$-cube (see
\cite{cubicalhdts}). This was the main motivation for introducing
cubical transition systems. Conversely, by
\cite[Proposition~9.7]{hdts}, there exists a labelled precubical set
$K$ such that its realization $\mathbb{T}(K)$ as weak transition
system does not satisfy CSA2. Every labelled precubical set is a
colimit of cubes, therefore $\mathbb{T}(K)$ is a colimit of cubes
since the realization functor from labelled symmetric precubical sets
to weak transition systems is colimit-preserving. Hence the weak
transition system $\mathbb{T}(K)$ is an example of a colimit of cubes
which is not regular (but it is cubical as any colimit of cubes).

\section{Intermediate state axiom and $\omega$-final lifts}
\label{section-cubical-lift}

Let $S$ be a set of objects of a locally presentable category $\K$.
For each object $X$ of $\K$, the colimit of the natural forgetful
functor $\widehat{S}\ddownarrow X \to \K$, where $\widehat{S}$ is the
full small category of $\K$ generated by $S$, is denoted by ($s\in S$
may be omitted) \[\liminj_{\begin{array}{c}s\rightarrow X\\s\in
    \mathcal{S}\end{array}}s.\] By
\cite[Proposition~3.1(i)]{1185.18014}, the full subcategory of
colimits of objects of $\mathcal{S}$ is a coreflective subcategory
$\K_\mathcal{S}$ of $\K$.  The right adjoint to the inclusion functor
$\K_\mathcal{S} \subset \K$ is precisely given by the functorial
mapping
\[X  \mapsto \liminj_{\begin{array}{c}s\rightarrow X\\s\in
      \mathcal{S}\end{array}}s.\]

  By \cite[Theorem~3.11]{cubicalhdts}, a weak transition system is
  cubical if and only if it is canonically a colimit of cubes and
  double transitions. In other terms, a weak transition system $X$ is
  cubical if and only if the canonical map
\[q_X : \liminj_{{\begin{array}{c}f: C_n[x_1,\dots,x_n] \rightarrow
      X\\f:\dd{x} \rightarrow X\end{array}}} \dom(f) \rightarrow X\]
is an isomorphism. The functorial mapping $X\mapsto \dom(q_X)$ is the
coreflection of the inclusion $\cts\subset \wts$.  The image of
$\underline{x}$ for any $x\in \Sigma$ by the coreflection $\wts \to
\cts$ is therefore the initial cubical transition system
$\varnothing$. This implies that the category $\cts$ is not a
concretely coreflective subcategory of $\wts$ over $\omega$ because
the set of actions is not preserved. Hence there is no reason for an
$\omega$-final lift of cubical transition systems to be cubical. This
holds anyway in the situation of Theorem~\ref{cubical-lift} which will
be used several times in the paper.

\bp \label{colimit-weak-CSA2} Let $X=\liminj X_i$ be a colimit of weak
transition systems. If all $X_i$ satisfy the Intermediate state axiom,
then so does $X$.  \ep

\bpf Let $T_i$ be the image by the canonical map $X_i\to X$ of the set
of transitions of $X_i$. Let $G_0=\bigcup_i T_i$. Let us define
$G_\alpha$ by induction on the transfinite ordinal $\alpha\geq 0$. If
$\alpha$ is a limit ordinal, then let $G_\alpha =
\bigcup_{\beta<\alpha}G_\beta$. If the set of tuples $G_\alpha$ is
defined, then let $G_{\alpha+1}$ be obtained from $G_\alpha$ by adding
the set of all $(q+2)$-tuples $(\nu_1,u_{p+1},\dots,u_{p+q},\nu_2)$
such that there exist five tuples $(\alpha,u_1,\dots,u_n,\beta)$,
$(\alpha,u_1,\dots,u_p,\nu_1)$, $(\nu_1,u_{p+1},\dots,u_n,\beta)$,
$(\alpha,u_1,\dots,u_{p+q},\nu_2)$ and
$(\nu_2,u_{p+q+1},\dots,u_n,\beta)$ of the set $G_\alpha$ for some
$p,q\geq 1$. For cardinality reason, the transfinite sequence
stabilizes and by \cite[Proposition~3.5]{hdts}, there exists an
ordinal $\alpha_0$ such that $G_{\alpha_0}$ is the set of transitions
of $X$. Every transition of $G_0$ satisfies the Intermediate state
axiom since it is satisfies by all $X_i$. Suppose that all transitions
of $G_\alpha$ satisfies the Intermediate state axiom. Take a tuple
$(\nu_1,u_{p+1},\dots,u_{p+q},\nu_2)$ of $G_{\alpha+1}$ like
above. Suppose that $q\geq 2$ and let $q>r\geq 1$. There exists a
state $\nu_3$ of $X$ such that the tuples
$(\alpha,u_1,\dots,u_{p+r},\nu_3)$ $(\nu_3,u_{p+r+1},\dots,u_n,\beta)$
are two transitions of $G_\alpha$ since all transitions of $G_\alpha$
satisfy the Intermediate state axiom by induction hypothesis.  From
the five tuples $(\alpha,u_1,\dots,u_n,\beta)$,
$(\alpha,u_1,\dots,u_{p+r},\nu_3)$,
$(\nu_3,u_{p+r+1},\dots,u_n,\beta)$ $(\alpha,u_1,\dots,u_{p+q},\nu_2)$
and $(\nu_2,u_{p+q+1},\dots,u_n,\beta)$ of $G_{\alpha}$, one deduces
that the tuple $(\nu_3,u_{p+r+1},\dots,u_{p+q},\nu_2)$ belongs to
$G_{\alpha+1}$. From the five tuples $(\alpha,u_1,\dots,\linebreak[4]u_n,\beta)$,
$(\alpha,u_1,\dots,u_p,\nu_1)$, $(\nu_1,u_{p+1},\dots,u_n,\beta)$,
$(\alpha,u_1,\dots,u_{p+r},\nu_3)$ and
$(\nu_3,u_{p+r+1},\dots,u_n,\linebreak[4]\beta)$, one deduces that the tuple
$(\nu_1,u_{p+1},\dots,u_{p+r},\nu_3)$ belongs to $G_{\alpha+1}$. Hence
$G_{\alpha+1}$ satisfies the Intermediate state axiom.  One deduces
that $X$ satisfies the Intermediate state axiom.  \epf

\bp \label{density-weak-CSA2} Consider the following map, functorial
with respect to the weak transition system $X$:
\[r_X : \liminj_{{\begin{array}{c}f: C_n[x_1,\dots,x_n] \rightarrow
      X\\f:\underline{x} \rightarrow X\end{array}}} \dom(f)
\rightarrow X.\] The map $r_X$ is always bijective on states and
actions and one-to-one on transitions.  The map $r_X$ is an
isomorphism if and only if $X$ satisfies the Intermediate state axiom.
\ep

\bpf Let $\alpha$ be a state of $X$. Then there exists a map $C_0\to
X$ sending the unique state of $C_0$ to $\alpha$. Hence $r_X$ is onto
on states. Let $\alpha$ and $\beta$ be two states of $\dom(r_X)$ sent
to the same state $\gamma$ of $X$.  Then the diagram $\{\alpha\}
\leftarrow \{\gamma\} \rightarrow \{\beta\}$ is a subdiagram in the
colimit calculating $\dom(r_X)$. Hence $\alpha = \beta$ in
$\dom(r_X)$. So $r_X$ is bijective on states. Let $u$ be an action of
$X$.  Then there exists a map $\underline{\mu(u)}\to X$ sending the
action $\mu(u)$ to $u$. This implies that $r_X$ is onto on actions.
Let $u$ and $v$ be two actions of $\dom(r_X)$ sent to the same action
$w$ of $X$.  Then the diagram $\{\mu(u)\} \leftarrow \{\mu(w)\}
\rightarrow \{\mu(v)\}$ is a subdiagram in the colimit calculating
$\dom(r_X)$. Hence $u=v$ in $\dom(r_X)$ and $r_X$ is bijective on
actions. Hence by \cite[Proposition~4.4]{homotopyprecubical}, $r_X$ is
always one-to-one on transitions.

By Proposition~\ref{colimit-weak-CSA2}, the weak transition system
$\dom(r_X)$ satisfies the Intermediate state axiom.  Therefore, if
$r_X$ is an isomorphism, then $X$ satisfies the Intermediate state
axiom. Conversely, let us suppose that $X$ satisfies the Intermediate
state axiom. Let $(\alpha,u_1,\dots,u_n,\beta)$ be a transition of
$X$. This transition gives rise to a map of weak transition systems
$\phi:C_n[\mu(u_1),\dots,\mu(u_n)]^{ext} \to X$. Since $X$ satisfies
the Intermediate state axiom, it is injective with respect to the
inclusion $C_n[\mu(u_1),\dots,\mu(u_n)]^{ext} \subset
C_n[\mu(u_1),\dots,\mu(u_n)]$ (see the proof of
\cite[Theorem~3.6]{cubicalhdts})~\footnote{Note that the composition
  axiom of weak transition systems is used here. It is worth noting
  that its use is often hidden.}. Hence $\phi$ factors as a composite
$C_n[\mu(u_1),\dots,\mu(u_n)]^{ext} \to C_n[\mu(u_1),\dots,\mu(u_n)]
\to X$. By definition of $\dom(r_X)$, $\phi$ factors as a composite
\[C_n[\mu(u_1),\dots,\mu(u_n)]^{ext} \longrightarrow
C_n[\mu(u_1),\dots,\mu(u_n)] \longrightarrow \dom(r_X) \stackrel{r_X}
\longrightarrow X.\] Hence $r_X$ is onto on transitions.  \epf

\bth \label{cubical-lift} Let $(f_i:\omega (A_i)\to W)_{i\in I}$ be a
cocone of $\set^{\{s\}\cup \Sigma}$ such that the weak transition
systems $A_i$ are cubical for all $i\in I$. Then the $\omega$-final
lift $\overline{W}$ satisfies the Intermediate state axiom.  Assume
moreover that every action $u$ of $W$ is the image of an action of
$A_{i_u}$ for some $i_u\in I$. Then the $\omega$-final lift
$\overline{W}$ is cubical.  \eth

\bpf Let $\C$ be the full subcategory of weak transition systems
satisfying the Intermediate axiom. By
Proposition~\ref{density-weak-CSA2} and
\cite[Proposition~3.1(i)]{1185.18014}, the category is a full
coreflective subcategory of $\wts$, the right adjoint being given by
the kelleyfication-like functor $X\mapsto \dom(r_X)$. Unlike the
coreflection from $\wts$ to $\cts$, the new coreflection preserves the
set of actions (and also the set of states). This means that the
category $\C$ is concretely coreflective over $\omega$. Hence
$\overline{W}$ satisfies the Intermediate state axiom by the dual of
\cite[Proposition~21.31]{topologicalcat}. Let $u$ be an action of
$\overline{W}$. Then, by hypothesis, there exists an action $v$ of
some $A_{i_u}$ such that the map $f_{i_u}:A_{i_u}\to W$ takes $v$ to
$u$. Since $A_{i_u}$ is cubical by hypothesis, there exists a
transition $(\alpha,v,\beta)$ of $A_{i_u}$. Hence the triple
$(f_{i_u}(\alpha),u,f_{i_u}(\beta))$ is a transition of
$\overline{W}$. This means that all actions of $\overline{W}$ are
used. In other terms, $\overline{W}$ is cubical.  \epf

Note that we have also proved that the forgetful functor $\C\subset
\wts\stackrel{\omega} \longrightarrow \set^{\{s\}\cup \Sigma}$ is
topological by \cite[Theorem~21.33]{topologicalcat}.  We give the
first application of this result. It states that the image of a
cubical transition system is cubical.

\begin{cor} \label{img-f} Let $f:X\to Y$ be a map of weak transition
  systems. Let $L_X$ ($L_Y$ resp.) be the set of actions of $X$ ($Y$
  resp.). Then $f$ factors as a composite $X\to f(X)\to Y$ such that
  the map $f(X)\to Y$ is the inclusion $f(X^0)\subset Y^0$ on states
  and the inclusion $f(L_X)\subset L_Y$ on actions. If $X$ is cubical,
  then $f(X)$ is cubical.  \end{cor}

\bpf Consider the $\omega$-final lift $f(X)$ of the map of
$\set^{\{s\}\cup \Sigma}$ \[\omega(X)\longrightarrow (f(X^0),f(L_X))\]
induced by $f$.  Then $f(X)$ is a solution. Assume now that $X$ is
cubical.  By Theorem~\ref{cubical-lift}, the weak transition system
$f(X)$ is cubical and the proof is complete. \epf

\section{Most elementary properties of regular transition systems}
\label{cat-rts}

A weak transition system satisfies the Unique intermediate state axiom
or CSA2 if and only if it is orthogonal to the set of inclusions
$C_n[x_1,\dots,x_n]^{ext} \subset C_n[x_1,\dots,x_n]$ for all $n\geq
0$ and all $x_1,\dots, x_n \in \Sigma$ by \cite[Theorem~5.6]{hdts}. By
\cite[Theorem~1.39]{MR95j:18001}, there exists a
functor \[\CSA_2:\wts\to \wts\] such that for any weak transition
system $Y$ satisfying CSA2 and any weak transition system $X$, the
weak transition system $\CSA_2(X)$ satisfies CSA2 and there is a
natural bijection $\wts(X,Y)\iso
\wts(\CSA_2(X),Y)$. Write \[\eta_X:X\to \CSA_2(X)\] for the unit of
this adjunction.  The following proposition provides an easy way to
check that a cubical transition system is regular.

\bp \label{mor-cub-reg} Let $X$ be a cubical transition system. Let
$Y$ be a weak transition system satisfying CSA2. Let $f:X\to Y$ be a
map of weak transition systems which is one-to-one on states.  Then
$X$ is regular.  \ep

Note that the hypothesis that $X$ is cubical cannot be
removed. Indeed, the inclusion \[C_n[x_1,\dots,x_n]^{ext} \subset
C_n[x_1,\dots,x_n]\] for $x_1,\dots, x_n \in \Sigma$ is one-to-one on
states because it is the inclusion $\{0_n,1_n\} \subset
\{0,1\}^n$. The target $C_n[x_1,\dots,x_n]$ satisfies CSA2. But the
pure $n$-transition $C_n[x_1,\dots,x_n]^{ext}$ does not satisfy CSA2
for $n\geq 2$ because it does not even satisfy the Intermediate state
axiom.

\bpf Let $(\alpha,u_1,\dots,u_n,\beta)$ be a transition of $X$ with
$n\geq 2$.  Let $1\leq p \leq n-1$. Since $X$ is cubical, there exist
two states $\nu_1$ and $\nu_2$ such that
$(\alpha,u_1,\dots,u_p,\nu_i)$ and $(\nu_i,u_{p+1},\dots,u_n,\beta)$
are transitions of $X$ for $i=1,2$. Then the five tuples
\begin{multline*}(f(\alpha),f(u_1),\dots,f(u_n),f(\beta)),\\
(f(\alpha),f(u_1),\dots,f(u_p),f(\nu_1)),
(f(\nu_1),f(u_{p+1}),\dots,f(u_n),f(\beta))\\
(f(\alpha),f(u_1),\dots,f(u_p),f(\nu_2)),
(f(\nu_2),f(u_{p+1}),\dots,f(u_n),f(\beta))
\end{multline*}
are transitions of $Y$.  Since $Y$ satisfies CSA2 by hypothesis, one
has $f(\nu_1)=f(\nu_2)$. Since $f$ is one-to-one on states by
hypothesis, one obtains $\nu_1=\nu_2$. Therefore $X$ satisfies CSA2.
\epf

\bp\label{calcul_CSA2} Let $X$ be a cubical transition system. There
exists a pushout diagram of cubical transition systems
\[
\xymatrix
{
X^0 \fR{\subset} \fD{(\eta_X)^0}&& X \fD{\eta_X} \\
&&\\
\CSA_2(X)^0 \fR{\subset} && \cocartesien \CSA_2(X)
}
\]
where the horizontal maps are the inclusion of the set of states into
the corresponding cubical transition system. For all cubical
transition systems $X$, the unit map $\eta_X:X\to \CSA_2(X)$ is onto
on states and the identity on actions.  \ep

Once again, the hypothesis that $X$ is cubical cannot be removed.
Indeed, let us consider again the case of a pure $n$-transition
$X=C_n[x_1,\dots,x_n]^{ext}$ with $x_1,\dots, x_n \in \Sigma$. Then
$\CSA_2(X)=C_n[x_1,\dots,x_n]$ by \cite[Theorem~5.6]{hdts}: in plain
English, the $n$-cube is the free regular transition system generated
the pure transition consisting of its $n!$ $n$-dimensional
transitions. The commutative square
\[
\xymatrix
{
\{0_n,1_n\} \fR{\subset} \fD{(\eta_{C_n[x_1,\dots,x_n]^{ext}})^0}&& C_n[x_1,\dots,x_n]^{ext} \fD{\eta_{C_n[x_1,\dots,x_n]^{ext}}} \\
&&\\
\{0,1\}^n \fR{\subset} && \cocartesien C_n[x_1,\dots,x_n]
}
\]
is not a pushout diagram. The unit map
$\eta_{C_n[x_1,\dots,x_n]^{ext}}$ is not onto on states. However, it
is still bijective on actions. 

We could actually prove that the map $\eta_X:X\to \CSA_2(X)$ is always
bijective on actions for any weak transition system $X$. We leave the
proof of this fact to the interested reader because it will not be
used in this paper.

\bpf The natural transformation from the state set functor
$(-)^0:\cts\to\set\subset \cts$ to the identity functor of $\cts$
gives rise to a commutative diagram of cubical transition systems:
\[
\xymatrix
{
X^0 \fR{\subset} \fD{(\eta_X)^0}&& X \fD{\eta_X} \\
&&\\
\CSA_2(X)^0 \fR{\subset} && \CSA_2(X).
}
\]
Consider the pushout diagram of cubical transition systems 
\[
\xymatrix
{
X^0 \fR{\subset} \fD{(\eta_X)^0}&& X \fD{\eta_X} \\
&&\\
\CSA_2(X)^0 \fR{} && \cocartesien Z.
}
\]
By the universal property of the pushout, the unit map $\eta_X:X\to
\CSA_2(X)$ factors uniquely as a composite \[X\longrightarrow
Z\longrightarrow \CSA_2(X).\] Since the forgetful functor $\omega:\wts
\to \set^{\{s\}\cup \Sigma}$ forgetting the transitions is
topological, and since the inclusion $\cts\subset \wts$ is
colimit-preserving, the state set functor $X\mapsto X^0$ from $\cts$
to $\set$ is colimit-preserving. Hence the set map $Z^0\to
\CSA_2(X)^0$ is bijective. Therefore, by
Proposition~\ref{mor-cub-reg}, the cubical transition system $Z$
satisfies CSA2. Hence we obtain $Z=\CSA_2(X)$ by the universal
property of the adjunction.  The functor taking a cubical transition
system to its set of actions is the composite functor
\[\xymatrix {
\cts \subset \wts \stackrel{\omega}\longrightarrow  \set^{\{s\}\cup \Sigma} \longrightarrow  \set^{\Sigma}\fR{L\mapsto 
\coprod_{x\in \Sigma}L_x} && \set}\] 
which is colimit-preserving as well. Therefore, one obtains the
pushout diagram of sets
\[
\xymatrix
{
\varnothing \fR{} \fD{}&& \hbox{set of actions of $X$} \fD{} \\
&&\\
\varnothing \fR{} && \cocartesien \hbox{set of actions of $\CSA_2(X)$}.
}
\]
This means that $X\to\CSA_2(X)$ is the identity on actions.  By
Corollary~\ref{img-f}, there exists a cubical transition system
$\eta_X(X)$ such that $\eta_X:X\to \CSA_2(X)$ factors as a composite
$X\to \eta_X(X)\to \CSA_2(X)$ such that the map $\eta_X(X)\to
\CSA_2(X)$ is the inclusion $\eta_X(X^0)\subset \CSA_2(X)^0$ on states
and an inclusion on actions.  By Proposition~\ref{mor-cub-reg},
$\eta_X(X)$ satisfies CSA2.  Therefore $\eta_X(X)=\CSA_2(X)$ by the
universal property of the adjunction. Hence the map $\eta_X:X\to
\CSA_2(X)$ is onto on states.  \epf

\bp \label{rts-csa2} If $X$ is cubical, then $\CSA_2(X)$ is
regular. In particular, if $X$ is regular, then $\CSA_2(X)$ is
regular. \ep

\bpf By definition, $\CSA_2(X)$ satisfies the Unique Intermediate
State axiom. By Proposition~\ref{calcul_CSA2}, the unit $X\to
\CSA_2(X)$ is the identity on actions. Therefore all actions of
$\CSA_2(X)$ are used since they are used in $X$ which is cubical.
\epf

\bp \label{rts-reflective} The category $\rts$ is a full reflective
subcategory of $\cts$ and the reflection is the functor
$\CSA_2:\cts\to\rts$ which is the restriction of $\CSA_2:\wts \to
\wts$ to cubical transition systems. \ep

\bpf Let $X$ be a cubical transition system and $Y$ a regular
transition system. By Proposition~\ref{rts-csa2}, one has the
bijection of sets $\cts(X,Y) \iso \rts(\CSA_2(X),Y)$. It is therefore
the left adjoint of the inclusion $\rts\subset \cts$.  \epf

\bp The category $\rts$ is locally finitely presentable. \ep

\bpf We already know that the cubes together with the double
transitions are a dense generator of $\cts$ by \cite[Theorem~3.11 and
Corollary~3.12]{cubicalhdts}. But they are regular. So $\rts$ has a
dense and hence strong generator because colimits in $\rts$ are
calculated, first, by taking the colimits in $\cts$ and, then, the
image by the reflection $\CSA_2:\cts\to\rts$. The category $\rts$ is
also cocomplete for the same reason. The proof is complete with
\cite[Theorem~1.20]{MR95j:18001}.  \epf

We can now introduce the cubification functor.

\bd \label{def-cub} \cite{hdts} \cite[Definition~3.13]{cubicalhdts}
Let $X\in \wts$. The {\rm cubification functor} is the functor $\cub :
\wts \longrightarrow \wts$ defined by \[\cub(X) =
\liminj_{C_n[x_1,\dots,x_n] \to X} C_n[x_1,\dots,x_n],\] the colimit
being taken in $\wts$.  \ed

For any $X\in \wts$, the weak transition system $\cub(X)$ is cubical
and the colimit can be taken in $\cts$ since the latter is
coreflective in $\wts$.

\bp \label{bij-state} Let $X$ be a weak transition system. Then the
canonical map \[\pi_X:\cub(X)\longrightarrow X\] is bijective on states. \ep

\bpf The argument is given in the proof of
\cite[Theorem~3.11]{cubicalhdts}.  \epf

\bp \label{cub-regular} Let $X$ be a regular transition system. Then
the cubical transition system $\cub(X)$ is regular and the colimit
\[\liminj_{C_n[x_1,\dots,x_n] \to X} C_n[x_1,\dots,x_n]\] is the same in
$\rts$, in $\cts$ and in $\wts$. \ep

\bpf The weak transition system $\cub(X)$ is cubical because it is a
colimit of cubes. The canonical map $\pi_X:\cub(X)\to X$ is bijective
on states by Proposition~\ref{bij-state}. Therefore $\cub(X)$ is
regular by Proposition~\ref{mor-cub-reg}. We already know that the
colimit is the same in $\cts$ and in $\wts$ since $\cts$ is a full
coreflective subcategory of $\wts$.  The functor $\CSA_2:\cts \to
\rts$ is a left adjoint to the inclusion $\rts\subset \cts$ by
Proposition~\ref{rts-reflective}. So it is colimit-preserving and one
obtains, because the cubes are regular, the isomorphism:
\[\CSA_2\left(\liminj\nolimits^{\cts} C_n[x_1,\dots,x_n]\right) \iso \liminj\nolimits^{\rts} C_n[x_1,\dots,x_n].
\]
The left-hand term is $\CSA_2(\cub(X))$ which is isomorphic to
$\cub(X)$ since $\cub(X)$ is regular.  \epf

\section{The left determined model category of regular transition systems}
\label{homotopy-rts}

Let us start this section with a few remarks about the terminology.

\begin{nota} For every map $f:X \to Y$ and every natural
  transformation $\alpha : F \to F'$ between two endofunctors of $\K$,
  the map $f\star \alpha$ is defined by the diagram:
\[
\xymatrix{
FX \fD{\alpha_X}\fR{f} && FY \fD{}\ar@/^15pt/@{->}[dddr]_-{\alpha_Y} &\\
&& &&\\
F'X \ar@/_15pt/@{->}[rrrd]^-{F'f}\fR{} &&  \bullet \cocartesien\ar@{->}[rd]^-{f\star \alpha} & \\
&& & F'Y.
}
\]
For a set of morphisms $\mathcal{A}$, let $\mathcal{A} \star \alpha =
\{f\star \alpha, f\in \mathcal{A}\}$.
\end{nota}

Let $(\C,\W,\F)$ be a model structure on a locally presentable
category $\K$ where $\C$ is the class of cofibrations, $\W$ the class
of weak equivalences and $\F$ the class of fibrations. A cylinder for
$(\C,\W,\F)$ is a triple $(\cyl:\K \to \K,\gamma^0\oplus \gamma^1:\id
\oplus \id\Rightarrow \cyl,\sigma:\cyl \Rightarrow \id)$ consisting of
a functor $\cyl:\K \to \K$ and two natural transformations
$\gamma^0\oplus \gamma^1:\id \oplus \id\Rightarrow \cyl$ and
$\sigma:\cyl \Rightarrow \id$ such that the composite $\sigma \circ
(\gamma^0\oplus \gamma^1)$ is the codiagonal functor $\id \oplus \id
\Rightarrow \id$ and such that the functorial map $\sigma_X:\cyl(X)
\to X$ belongs to $\W$ for every object $X$.  We will often use the
notation $\gamma=\gamma^0\oplus \gamma^1$.  The cylinder is
\emph{good} if the functorial map $\gamma_X : X \sqcup X \to \cyl(X)$
is a cofibration for every object $X$. It is \emph{very good} if,
moreover, the map $\sigma_X:\cyl(X) \to X$ is a trivial fibration for
every object $X$. A good cylinder is \emph{cartesian} if
\begin{itemize}
\item The functor $\cyl:\K\to \K$ has a right adjoint $\cocyl:\K \to
  \K$ called the {\rm path functor}.
\item There are the inclusions $\C \star \gamma^\epsilon \subset \C$
  for $\epsilon=0,1$ and $\C \star \gamma \subset \C$.
\end{itemize}
The notions above can be adapted to a cofibrantly generated weak
factorization system $(\mathcal{L},\mathcal{R})$ by considering the
combinatorial model structure
$(\mathcal{L},\Mor(\K),\mathcal{R})$. They can be also extended to any
set of maps $I$ by considering the associated cofibrantly generated
weak factorization system in the sense of
\cite[Proposition~1.3]{MR1780498}.

\bd \label{boundary-def} Let $n\geq 1$ and $x_1,\dots,x_n \in
\Sigma$. Let $\de C_n[x_1,\dots,x_n]$ be the regular transition system
defined by removing from the $n$-cube $C_n[x_1,\dots,x_n]$ all its
$n$-transitions. It is called the {\rm boundary} of
$C_n[x_1,\dots,x_n]$. \ed

\begin{nota} \label{cofgen} Denote by $\I$ the set of maps of cubical
  transition systems:
\begin{multline*}
\I = \{C:\varnothing \to \{0\}, R:\{0,1\} \to \{0\}\} \\ \cup \{\de C_n[x_1,\dots,x_n]
\to C_n[x_1,\dots,x_n]\mid \hbox{$n\geq 1$ and $x_1,\dots,x_n \in
  \Sigma$}\} \\ \cup \{C_1[x] \to \dd{x}\mid x\in \Sigma\}.
\end{multline*}
\end{nota}

By \cite[Corollary~6.8]{cubicalhdts} and
\cite[Theorem~4.6]{homotopyprecubical}, there exists a (necessarily
unique) left determined model category structure on $\cts$ (denoted by
$\cts$ as well) with the set of generating cofibrations $\I$. A map of
cubical transition systems is a cofibration of this model structure if
and only if it is one-to-one on actions. By
\cite[Proposition~5.5]{cubicalhdts}, this model category has a
cartesian and very good cylinder $\cyl:\cts\to \cts$ defined on
objects as follows: for a cubical transition system $X=(S,\mu:L\to
\Sigma,T)$, $\cyl(X)$ has the same set of states $S$, the set of
actions $L\p \{0,1\}$ with the labelling map $L\p \{0,1\} \to L \to
\Sigma$ and a tuple $(\alpha,(u_1,\epsilon_1), \dots,
(u_n,\epsilon_n), \beta)$ is a transition of $\cyl(X)$ if and only if
$(\alpha,u_1, \dots,u_n, \beta)$ is a transition of $X$. The map
$\gamma^\epsilon_X:X\to \cyl(X)$ for $\epsilon=0,1$ is induced by the
identity on states and by the mapping $u\mapsto (u,\epsilon)$ on
actions. The map $\sigma_X:\cyl(X)\to X$ is induced by the identity on
states and by the projection $(u,\epsilon)\mapsto u$ on actions.

\bp \label{cyl-rts} One has the natural isomorphism of cubical
transition systems
\[\CSA_2(\cyl(X)) \iso \cyl(\CSA_2(X))\] for every cubical transition
system $X$. \ep

\bpf We have just recalled that the canonical map $\sigma_X:\cyl(X)\to
X$ is bijective on states.  Therefore, by
Proposition~\ref{mor-cub-reg}, one has $\cyl(\rts)\subset \rts$. By
Proposition~\ref{calcul_CSA2}, for every cubical transition system
$X$, one has the pushout diagram of weak transition systems (and of
cubical transition systems since colimits are the same):
\[
\xymatrix
{
X^0 \fR{}\fD{} && X \fD{} \\
&&\\
\CSA_2(X)^0 \fR{} && \cocartesien \CSA_2(X).}
\]
Since $\cyl:\cts\to \cts$ is a left adjoint, one obtains the pushout
diagram of cubical transition systems:
\[
\xymatrix
{
\cyl(X^0) \fR{}\fD{} && \cyl(X) \fD{} \\
&&\\
\cyl(\CSA_2(X)^0) \fR{} && \cocartesien \cyl(\CSA_2(X)).}
\]
For any set $E$ viewed as a cubical transition system, one has
$\cyl(E)=E$. Therefore one obtains the pushout diagram of cubical
transition systems:
\[
\xymatrix
{
X^0 \fR{}\fD{} && \cyl(X) \fD{} \\
&&\\
\CSA_2(X)^0 \fR{} && \cocartesien \cyl(\CSA_2(X)).}
\]
Since $\CSA_2(X)$ is regular, the cubical transition system
$\cyl(\CSA_2(X))$ is regular. Therefore, by
Proposition~\ref{calcul_CSA2}, the cubical transition system
$\cyl(\CSA_2(X))$ and $\CSA_2(\cyl(X))$ satisfy the same universal
property. Hence we obtain the natural isomorphism \[\CSA_2(\cyl(X))
\iso \cyl(\CSA_2(X)).\] \epf

\bth \label{left-determined-rhdts} There exists a (necessarily unique)
left determined model category structure on $\rts$ (denoted by $\rts$)
such that the set of generating cofibrations is $\CSA_2(\I)=\I$ and
such that the fibrant objects are the fibrant cubical transition
systems which are regular. The cartesian cylinder is the restriction
to $\rts$ of the cylinder of $\cts$ defined above. The restricted
cylinder is very good. The reflection $\CSA_2:\cts \to \rts$ is a left
Quillen homotopically surjective functor. The inclusion $\rts \subset
\cts$ reflects weak equivalences. \eth

\bpf Thanks to Proposition~\ref{reflective-adjunction} applied with
Proposition~\ref{cyl-rts}, we see that $\cyl:\cts\to\cts$ and its
right adjoint $\cocyl:\cts\to\cts$ restrict to endofunctors of
$\rts$. We then apply \cite[Lemma~5.2]{MO} which is reexplained also
in \cite[Theorem~9.3]{homotopyprecubical}. The only thing which
remains to be proved is that the restriction $\cyl:\rts\to\rts$ is a
very good cylinder.  Consider the following commutative square of
solid arrows of $\rts$:
\[
\xymatrix{
\CSA_2(A)\fR{}\fD{\CSA_2(f)} && \cyl(X) \fD{\sigma_X}\\
&&\\
\CSA_2(B) \fR{}\ar@{-->}[rruu]^-{k}&& X}
\]
where $f\in \I$ and $X\in\rts$. Because of the adjunction, the existence
of a lift $k$ is equivalent to the existence of a lift in the
following commutative square of solid arrows of $\cts$:
\[
\xymatrix{
A\fR{}\fD{f} && \cyl(X) \fD{\sigma_X}\\
&&\\
B \fR{}\ar@{-->}[rruu]^-{\ell}&& X.}
\]
So the restriction of $\cyl$ to $\rts$ is very good as well.
\epf

The end of the section is devoted to a characterization of the weak
equivalences of the left-determined model structure $\rts$.

\bp \label{rts-fibrant} (Compare with
\cite[Proposition~7.8]{cubicalhdts}) Every regular transition system
satisfying CSA1 is fibrant in $\rts$. The category of regular
transition systems satisfying CSA1 is a small-orthogonality class of
$\rts$. \ep

\bpf Every regular transition system satisfying CSA1 is fibrant in
$\cts$ by \cite[Proposition~7.8]{cubicalhdts}, and therefore fibrant
in $\rts$ by Corollary~\ref{left-determined-rhdts}. A regular
transition system is CSA1 if and only if it is orthogonal to the maps
$\sigma_{C_1[x]}:\cyl(C_1[x])\to C_1[x]$ for all $x\in \Sigma$.  \epf

The full subcategory of regular transition systems satisfying CSA1 is
therefore a full reflective subcategory by
\cite[Theorem~1.39]{MR95j:18001}. Write
$\CSA_1^{\rts}:\rts\longrightarrow\rts$ for the reflection.  The full subcategory
of cubical transition systems satisfying CSA1 is also a
small-orthogonality class and a full reflective subcategory of $\cts$
by
\cite[Proposition~7.2]{cubicalhdts}. Write $\CSA_1^{\cts}:\cts\longrightarrow\cts$
for the reflection. The functor $\CSA_1^{\rts}:\rts\to \rts$
($\CSA_1^{\cts}:\cts\to\cts$ resp.)  can be defined as follows. Let
$X_0=X$. We construct by transfinite induction a sequence of regular
(cubical resp.) transition systems as follows: if for $\alpha\geq 0$,
there exist two transitions $(\alpha,u,\beta)$ and $(\alpha,u',\beta)$
with $u\neq u'$ and $\mu(u)=\mu(u')$, consider the pushout diagram in
$\rts$ (in $\cts$ resp.)
\[
\xymatrix
{
\cyl(C_1[\mu(u)]) \fR{{\tiny\begin{array}{c}(\mu(u),1,0)\mapsto u\\ (\mu(u),1,1)\mapsto u'\end{array}}} \fD{\sigma_{C_1[\mu(u)]}} && X_\alpha \fD{} \\
&& \\
C_1[\mu(u)] \fR{} && \cocartesien X_{\alpha+1},
}
\]
otherwise let $X_{\alpha+1}=X_\alpha$. If $\alpha$ is a limit ordinal,
then let $X_\alpha=\liminj_{\beta<\alpha} X_\beta$, the colimit being
calculated $\rts$ (in $\cts$ resp.). By a cardinality argument
(all maps $X_\alpha\to X_{\alpha+1}$ are onto on actions), the
sequence stabilizes.  The colimit is $\CSA_1^{\rts}(X)$
($\CSA_1^{\cts}(X)$ resp.). 

Let $X$ be a regular transition system. The canonical map $X\to
\CSA_1^{\cts}(X)$ is then a transfinite composition of pushouts in
$\cts$ of maps of $\{\sigma_{C_1[x]}\mid x\in \Sigma\}$. Since a
colimit is calculated in $\rts$ by taking the colimit in $\cts$ and by
taking the image by the functor $\CSA_2$, the map $\CSA_2(X)=X\to
\CSA_2(\CSA_1^{\cts}(X))$ is a transfinite composition of pushouts in
$\rts$ of maps of $\{\sigma_{C_1[x]}\mid x\in \Sigma\}$. Thus, 
$\CSA_1^{\rts}(X)$ is orthogonal to $\CSA_2(X)=X\to
\CSA_2(\CSA_1^{\cts}(X))$. Hence the canonical map $X\to
\CSA_1^{\rts}(X)$ factors uniquely as a composite \[X\longrightarrow
\CSA_2(\CSA_1^{\cts}(X)) \longrightarrow \CSA_1^{\rts}(X).\]

\bp \label{CSA1-diff} There exists a regular transition system $X$
such that the ``comparison map''
\[\CSA_2(\CSA_1^{\cts}(X)) \to \CSA_1^{\rts}(X)\] is not an
isomorphism. \ep 

\bpf A cubical transition system is completely defined by giving the
list of all transitions and the actions identified by the labelling
map.  We consider the regular transition system $X$ having the
transitions
\begin{multline*}
(\alpha,u_1,u_2,\beta),(\alpha,u_2,u_1,\beta),(\alpha,u_1,\chi),(\chi,u_2,\beta),(\alpha,u_2,\nu),(\nu,u_1,\beta),\\
(\alpha,u'_1,u'_2,\beta),(\alpha,u'_2,u'_1,\beta),(\alpha,u'_1,\chi'),(\chi',u'_2,\beta),(\alpha,u'_2,\nu'),(\nu',u'_1,\beta),\\
(\gamma,v,\chi),(\gamma,v',\chi'),(U_1,u_1,V_1),(U_1,u'_1,V_1),(U_2,u_2,V_2),(U_2,u'_2,V_2)
\end{multline*}
such that all actions are labelled by some $x\in\Sigma$. By applying
the functor $\CSA_1^{\cts}:\cts \to \cts$ to $X$, the actions $u_i$
and $u'_i$ are identified because of the presence of the transitions
\[(U_1,u_1,V_1),(U_1,u'_1,V_1),(U_2,u_2,V_2),(U_2,u'_2,V_2).\] The
functor $\CSA_1^{\cts}:\cts \to \cts$ does not make the identification
$v=v'$ because these two actions are used in the transitions
$(\gamma,v,\chi)$ and $(\gamma,v',\chi')$ and because it is assumed
that $\chi\neq \chi'$.  The cubical transition system
$\CSA_1^{\cts}(X)$ therefore consists of the transitions~\footnote{The
  states are preserved by $\CSA_1^{\cts}$ since the canonical map
  $X\to\CSA_1^{\cts}(X)$ is a transfinite composition of pushouts of
  maps of the form $\cyl(C_1[z])\to C_1[z]$ for $z\in \Sigma$, because
  these maps are all of them state-preserving and because the state
  set functor from $\cts$ to $\set$ is colimit-preserving. Beware of
  the fact that the functor $\CSA_1^{\rts}$ is not state-preserving.}
\begin{multline*}
(\alpha,u_1,u_2,\beta),(\alpha,u_2,u_1,\beta),(\alpha,u_1,\chi),(\chi,u_2,\beta),(\alpha,u_2,\nu),(\nu,u_1,\beta),\\
(\alpha,u_1,u_2,\beta),(\alpha,u_2,u_1,\beta),(\alpha,u_1,\chi'),(\chi',u_2,\beta),(\alpha,u_2,\nu'),(\nu',u_1,\beta),\\
(\gamma,v,\chi),(\gamma,v',\chi'),(U_1,u_1,V_1),(U_2,u_2,V_2).
\end{multline*}
The latter cubical transition system is not regular. Indeed, in the
regular transition system $\CSA_2(\CSA_1^{\cts}(X))$, the
identifications of states $\chi=\chi'$ and $\nu=\nu'$ are made. We
obtain for $\CSA_2(\CSA_1^{\cts}(X))$ the list of transitions
\begin{multline*}
(\alpha,u_1,u_2,\beta),(\alpha,u_2,u_1,\beta),(\alpha,u_1,\chi),(\chi,u_2,\beta),(\alpha,u_2,\nu),(\nu,u_1,\beta),\\
(\gamma,v,\chi),(\gamma,v',\chi),(U_1,u_1,V_1),(U_2,u_2,V_2).
\end{multline*}
The map $\CSA_2(\CSA_1^{\cts}(X)) \to \CSA_1^{\rts}(X)$ therefore
identifies the actions $v$ and $v'$. Hence it is not an isomorphism.
\epf

\bp \label{homotopy_is_equality} (Compare with
\cite[Proposition~7.4]{cubicalhdts}) Let $Y$ be a regular transition
system satisfying CSA1.  Let $X$ be a regular transition system. Then
two homotopy equivalent maps $f,g:X\rightarrow Y$ are equal. In other
terms, each of the two canonical maps $X \rightarrow \cyl(X)$ induces a
bijection $\rts(\cyl(X),Y) \iso \rts(X,Y)$. \ep

\bpf By \cite[Proposition~7.4]{cubicalhdts}, one has the bijection of
sets \[\cts(\cyl(X),Y) \iso \cts(X,Y),\] the binary product being
calculated in $\cts$.  The category $\rts$ is a full reflective
subcategory of $\cts$ by Proposition~\ref{rts-reflective}. Thus, there
is the bijection \[\rts(\cyl(X),Y) \iso \rts(X,Y)\] where the binary
product is calculated in $\rts$.  \epf

The following model-categorical lemma is implicitly used several
times in \cite{cubicalhdts} and \cite{homotopyprecubical} and it will
be used again several times in this paper. Let us state it clearly:

\begin{lem} \label{when-cell-is-weak} Let $\mathcal{M}$ be a left
  proper combinatorial model category such that the generating
  cofibrations are maps between finitely presentable objects. Let $\C$
  be a class of weak equivalences of $\mathcal{M}$ satisfying the
  following condition: in every pushout diagram of $\mathcal{M}$ of
  the form
\[
\xymatrix
{
A \fD{g\in \C} \fR{\phi} && C \fD{f}\\
&& \\
B \fR{}&& \cocartesien D
}
\]
either $\phi$ is a cofibration or $f$ is an isomorphism.  Then every
map of $\cell_{\mathcal{M}}(\C)$ is a weak equivalence of
$\mathcal{M}$, where $\cell_{\mathcal{M}}(\C)$ is the class of
transfinite composition of pushouts of maps of $\C$.
\end{lem}

\bpf Since $\mathcal{M}$ is left proper, $f$ is always a weak
equivalence of $\mathcal{M}$.  By \cite[Proposition~4.1]{rankweak},
the class of weak equivalences of $\mathcal{M}$ is closed under
transfinite composition.  Hence the proof is complete. \epf

\begin{lem} \label{for-when-cell-is-weak-1} For all $x\in \Sigma$, the
  map $\sigma_{C_1[x]}:\cyl(C_1[x])\to C_1[x]$ satisfies the
  conditions of Lemma~\ref{when-cell-is-weak} for $\mathcal{M}=\rts$.
\end{lem}

\bpf Consider a pushout diagram of $\rts$ 
\[
\xymatrix
{
\cyl(C_1[x]) \fD{\sigma_{C_1[x]}} \fR{\phi} && C \fD{f}\\
&& \\
C_1[x] \fR{}&& \cocartesien D.
}
\]

The map $f:C\to D$ factors as a composite $f:C\to E\to \CSA_2(E)=D$
where $E$ is the colimit in $\cts$. If $\phi$ is not a cofibration,
then $\phi$ is constant on actions.  In this case, $C\iso E$ by the
proof of \cite[Theorem~7.10]{cubicalhdts}, therefore $E$ is regular.
One obtains $D=\CSA_2(E)\iso E\iso C$. Hence $f$ is an isomorphism.
\epf

\bth \label{weak-rts} (Compare with \cite[Theorem~7.10]{cubicalhdts})
A map $f:X\to Y$ of regular transition systems is a weak equivalence
for the left determined model structure of $\rts$ if and only if the
map $\CSA^{\rts}_1(f):\CSA^{\rts}_1(X) \to \CSA^{\rts}_1(Y)$ is an
isomorphism. \eth

\bpf By Lemma~\ref{for-when-cell-is-weak-1}, a map of regular
transition systems $f:X\to Y$ is a weak equivalence if and only if the
map $\CSA^{\rts}_1(f):\CSA^{\rts}_1(X) \to \CSA^{\rts}_1(Y)$ is a weak
equivalence. Since $\CSA^{\rts}_1(X)$ and $\CSA^{\rts}_1(Y)$ are
fibrant by Proposition~\ref{rts-fibrant}, a map of regular transition
systems $f:X\to Y$ is a weak equivalence if and only if the map
$\CSA^{\rts}_1(f):\CSA^{\rts}_1(X) \to \CSA^{\rts}_1(Y)$ is a homotopy
equivalence. The proof is complete with
Proposition~\ref{homotopy_is_equality}.  \epf

\section{Bousfield localization of the regular t.s.  by the
  cubification functor}
\label{Bousfield-cub}

We now deal with the Bousfield localization of $\rts$ by the
cubification functor $\cub$ and we compare this Bousfield localization
with the one of $\cts$ by the same cubification functor.

Let $x\in \Sigma$. Consider the unique map $p_x:C_1[x] \sqcup
C_1[x]\to \dd{x}$ bijective on states and sending the actions of the
source $C_1[x] \sqcup C_1[x]$ to their label.  Let us factor $p_x$ as
a composite (all maps are bijective on states)
\[
\xymatrix{C_1[x] \sqcup
  C_1[x]\ar@{^(->}[rr]_-{p_x^{cof}}^-{{\tiny\begin{array}{c}x_1\mapsto
        x_1\\x_2\mapsto x_2\end{array}}} && Z^{x_1,x_2}_x
  \ar@{->}[rr]_-{\simeq}^-{{\tiny\begin{array}{c}x_1\mapsto
        x_1\\x_2\mapsto x_2\\x\mapsto x\end{array}}} && \dd{x}}\] with
$Z^{x_1,x_2}_x$ is depicted in Figure~\ref{Z^{x_1,x_2}_x}, and where
$x_1$ and $x_2$ are the two actions of $C_1[x] \sqcup C_1[x]$ with
$\mu(x_1)=\mu(x_2)=x$. The left-hand map is a cofibration because it
is one-to-one on actions. One has the
isomorphisms \[\CSA^{\rts}_1(Z^{x_1,x_2}_x) \iso \CSA^{\rts}_1(\dd{x})
\iso\CSA^{\cts}_1(Z^{x_1,x_2}_x) \iso \CSA^{\cts}_1(\dd{x}) \iso
\dd{x},\] so the right-hand map is a weak equivalence of $\cts$ by
\cite[Theorem~7.10]{cubicalhdts} and of $\rts$ by
Theorem~\ref{weak-rts}. Therefore $p_x^{cof}$ is a cofibrant
replacement of $p_x$ both in $\cts$ and in $\rts$.

\begin{nota} Let $\mathcal{S}=\{p_x\mid x\in \Sigma\}$ and
  $\mathcal{S}^{cof}=\{p_x^{cof}\mid x\in \Sigma\}$.
\end{nota}

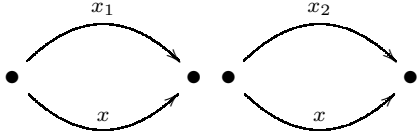
\begin{figure}
\[
\xymatrix{
\bullet \ar@/^20pt/[rr]^-{x_1} \ar@/_20pt/[rr]^-{x}&& \bullet }
\xymatrix{
\bullet \ar@/^20pt/[rr]^-{x_2} \ar@/_20pt/[rr]^-{x}&& \bullet }
\]
\caption{The cubical transition system $Z^{x_1,x_2}_x$ contains four states and three actions $x_1,x_2,x$ with $\mu(x_1)=\mu(x_2)=x$.}
\label{Z^{x_1,x_2}_x}
\end{figure}

\bp \label{simpl} For a cubical transition system $X$, the following
statements are equivalent:
\begin{enumerate}
\item The labelling map $\mu$ is one-to-one.
\item $X$ is $\mathcal{S}$-injective.
\item $X$ is $\mathcal{S}$-orthogonal.
\end{enumerate}
If any one of these statements is true, then $X$ satisfies CSA1 and is
$\mathcal{S}^{cof}$-orthogonal.  \ep

\bpf The equivalence $(1)\Longleftrightarrow (2) \Longleftrightarrow
(3)$ and the fact that these three conditions imply CSA1 is
\cite[Proposition~8.2]{homotopyprecubical}. Let $X$ be a cubical
transition system satisfying $(1)$. Consider the diagram of cubical
transition systems:
\[
\xymatrix
{
C_1[x]\sqcup C_1[x] \fD{p_x^{cof}} \fR{\phi} && X \\
&& \\
Z^{x_1,x_2}_x, \ar@{-->}[rruu]^-{\ell}&&
}
\]
where $x_1$ and $x_2$ are the two actions of $C_1[x]\sqcup C_1[x]$.
Define $\ell$ on states by $\ell(\alpha)=\phi(\alpha)$ for all states
$\alpha$, and on actions by $\ell(x_i)=\phi(x_i)$ for $i=1,2$ and
$\ell(x)=\phi(x_1)$.  Since $X$ satisfies $(1)$, one has
$\phi(x_1)=\phi(x_2)$. We deduce that $\ell$ is a well-defined map of
cubical transition systems. The map $\ell$ is the only solution
because $p_x^{cof}$ is bijective on states and the image by $\ell$ of
the actions of $Z^{x_1,x_2}_x$ is necessarily the unique action of $X$
labelled by $x$. Hence $X$ is $\mathcal{S}^{cof}$-orthogonal.  \epf

\bp \label{pi-cell-rts} (Compare with
\cite[Proposition~8.4]{cubicalhdts}) For every regular transition
system $X$, the canonical map $\pi_X:\cub(X)\to X$ belongs to
$\cell_{\rts}(\mathcal{S})$ \ep

\bpf The difficulty is, once again, that colimits are not calculated
in the same way in $\rts$ and in $\cts$. Let $(u_1^i,u_2^i)_{i\in I}$
be the family of pairs of actions of $X$ such that $\pi_X(u_1^i) =
\pi_X(u_2^i)$, which implies $\mu(u_1^i)=\mu(u_2^i)$. Since $X$ is
cubical, for all $i\in I$, there exist $1$-transitions
$(\alpha_j^i,u_j^i,\beta_j^i)$ of $X$ for $j=1,2$. Let $\phi^i:
C_1[\mu(u_1^i)] \sqcup C_1[\mu(u_2^i)]\to X$ be the map of cubical
transition systems sending the two $1$-transitions of the source to
$(\alpha_j^i,u_j^i,\beta_j^i)$ for $j=1,2$. Since $\pi_X:\cub(X)\to X$
is the identity on states by Proposition~\ref{bij-state}, one obtains
the following commutative diagram of regular transition systems:
\[
\xymatrix{
\coprod\limits_{i\in I} C_1[\mu(u_1^i)] \sqcup C_1[\mu(u_2^i)] \fR{\phi^i} \fD{\coprod\limits_{i\in I}p_{\mu(u_1^i)}}&& \cub(X) \fD{\pi_X} \\
&& \\
\coprod\limits_{i\in I}\dd{\mu(u_1^i)}\fR{} && X.}
\]
Consider the pushout diagram of regular transition systems:
\[
\xymatrix{
\coprod\limits_{i\in I} C_1[\mu(u_1^i)] \sqcup C_1[\mu(u_2^i)] \fR{\phi^i} \fD{\coprod\limits_{i\in I}p_{\mu(u_1^i)}}&& \cub(X) \fD{} \\
&& \\
\coprod\limits_{i\in I}\dd{\mu(u_1^i)}\fR{} && \cocartesien Z.}
\]
The colimit $Z$ is calculated in $\rts$ by taking the colimit $T$ in
$\cts$ and by taking the image by the reflection $\CSA_2$. Hence the
map $\pi_X:\cub(X)\to X$ factors as a
composite \[\cub(X)\longrightarrow T \longrightarrow \CSA_2(T)=Z
\stackrel{h}\longrightarrow X.\] The map $\cub(X)\to T$ is a pushout
in $\cts$ of the map $\coprod_{i\in I}p_{\mu(u_1^i)}$. The latter is
bijective on states, therefore the map $\cub(X)\to T$ is bijective on
states as well. The map $T\to Z$ is onto on states by
Proposition~\ref{calcul_CSA2}. Hence the map $g:\cub(X)\to Z$ is onto
on states. Let $\alpha$ and $\beta$ be two states of $\cub(X)$ mapped
to the same state $\gamma$ of $Z$. Then $\gamma$ is mapped to
$\pi_X(\alpha)=\pi_X(\beta)$ by $Z\to X$. Hence $\alpha=\beta$ by
Proposition~\ref{bij-state}. Therefore $g:\cub(X)\to Z$ is bijective
on states, and so is the map of cubical transition systems $h:Z\to
X$. By construction, the latter map is one-to-one on
actions. Therefore $h:Z\to X$ is one-to-one on transitions by
\cite[Proposition~4.4]{homotopyprecubical}.  Any action $u$ is used by
a $1$-transition $(\alpha,u,\beta)$ of $X$.  Hence $\pi_X:\cub(X)\to
X$ is onto on actions. Thus, there exists an action $v$ of $\cub(X)$
such that $\pi_X(v)=u$. This means that $h(g(v))=u$.  Hence $h$ is
onto on actions as well. To conclude that $h$ is an isomorphism,
consider a transition $(\alpha,u_1,\dots,u_n,\beta)$ of $X$. It gives
rise to a map of weak transition systems
$C_n[\mu(u_1),\dots,\mu(u_n)]^{ext}\to X$ which factors as a composite
$C_n[\mu(u_1),\dots,\mu(u_n)]^{ext}\to C_n[\mu(u_1),\dots,\mu(u_n)]
\to X$ since $X$ is cubical.  One obtains the composite map of weak
transition systems
\[C_n[\mu(u_1),\dots,\mu(u_n)]^{ext} \longrightarrow
C_n[\mu(u_1),\dots,\mu(u_n)] \longrightarrow \cub(X) \longrightarrow
Z\longrightarrow X.\] Hence $h$ is onto on transitions.
\epf

\begin{lem} \label{for-when-cell-is-weak-2} For all $x\in \Sigma$, the
  map $p_x:C_1[x]\sqcup C_1[x]\to \dd{x}$ satisfies the conditions of
  Lemma~\ref{when-cell-is-weak} for $\mathcal{M}=\bl_{\cub}\rts$.
\end{lem}

\bpf  Consider a pushout diagram of $\rts$ 
\[
\xymatrix
{
C_1[x]\sqcup C_1[x] \fD{p_x} \fR{\phi} && C \fD{f}\\
&& \\
\dd{x} \fR{}&& \cocartesien D.
}
\]
The map $f:C\to D$ factors as a composite $f:C\to E\to \CSA_2(E)=D$
where $E$ is the colimit in $\cts$. If $\phi$ is not a cofibration,
then $\phi$ is constant on actions.  In this case, $C\iso E$ by the
proof of \cite[Proposition~8.5]{cubicalhdts}, therefore $E$ is
regular.  One obtains $D=\CSA_2(E)\iso E\iso C$. Hence $f$ is an
isomorphism.  \epf

\bth (Compare with
\cite[Theorem~8.6]{cubicalhdts}) \label{eg-localizer} Let $\W_{\cub}$
be the Grothendieck localizer generated by the class of maps $f:X\to
Y$ of regular transition systems such that $\cub(f):\cub(X)\to
\cub(Y)$ is a weak equivalence of $\rts$ (the left determined model
structure). Let $\W(\mathcal{S})$ be the Grothendieck localizer
generated by the set of maps $\mathcal{S}$.  Then one has
$\W_{\cub}=\W(\mathcal{S})$. \eth

\bpf The proof is \emph{mutatis mutandis} the proof of
\cite[Theorem~8.6]{cubicalhdts}. Let us sketch it.  By
Proposition~\ref{pi-cell-rts}, the counit $\pi_X:\cub(X)\to X$ belongs
to $\cell_{\rts}(\mathcal{S})$ for all regular transition systems.  By
Lemma~\ref{for-when-cell-is-weak-2}, one deduces that
$\cell_{\rts}(\mathcal{S}) \subset \W(\mathcal{S})$. Hence, for all
regular transition systems $X$, the counit $\pi_X:\cub(X)\to X$
belongs to $\W(\mathcal{S})$.  Let $f:X\to Y$ be a map of
$\W_{\cub}$. Consider the commutative diagrams:
\[
\xymatrix{
\cub(X) \fR{\cub(f)} \fD{} && \cub(Y) \fD{} \\
&& \\
X \fR{f} && Y.}
\]
We have just proved that the vertical maps belong to
$\W(\mathcal{S})$.  Since $\cub(f)$ is a weak equivalence of $\rts$,
i.e. it belongs to the smallest Grothendieck localizer
$\W(\varnothing)\subset \W(\mathcal{S})$, one deduces by the
two-out-of-three property that the bottom map $f$ belongs to
$\W(\mathcal{S})$ as well. Hence we obtain the inclusion $\W_{\cub}
\subset \W(\mathcal{S})$.  Since $\cub(p_x)$ is an automorphism of
$C_1[x]\cup C_1[x]$, one has $\mathcal{S}\subset\W_{\cub}$, and
therefore $\W(\mathcal{S})\subset \W_{\cub}$.  \epf

\begin{cor} (Compare with
  \cite[Corollary~8.7]{cubicalhdts}) \label{localization} The
  Bousfield localization of the left determined model structure of
  $\rts$ with respect to the functor $\cub$ exists. \end{cor}

\bpf The combinatorial model category $\rts$ is left proper since all
objects are cofibrant. We want to Bousfield localize with respect to a
set of maps $\mathcal{S}$. Hence the proof is complete.  \epf

\begin{nota} Let us write $\bl_{\cub}\cts$ ($\bl_{\cub}\rts$ resp.)
  for the Bousfield localization of $\cts$ ($\rts$ resp.) by the
  functor $\cub$. \end{nota}

\bp \label{fibrant-rts} A regular transition system is fibrant in
$\bl_{\cub}\rts$ if and only if it is fibrant in $\bl_{\cub}\cts$. \ep

\bpf The proof is similar to the proof of
Theorem~\ref{left-determined-rhdts}.  \epf

\bp (Compare with \cite[Theorem~8.11 (1)(2)(3)]{cubicalhdts})
\label{def-LS-rts} The category
\[\inj_{\rts}(\mathcal{S})\] of $\mathcal{S}$-injective regular
transition systems is a small-orthogonality class and a full
reflective subcategory of $\rts$. Write $\bl^{\rts}_{\mathcal{S}} :
\rts \rightarrow \rts$ for the reflection. The unit map $X\to
\bl^{\rts}_{\mathcal{S}}(X)$ belongs to $\cell_{\rts}(\mathcal{S})$
for any regular transition system $X$.  \ep

\bpf By Proposition~\ref{simpl}, being $\mathcal{S}$-injective is
equivalent to being $\mathcal{S}$-orthogonal.  By
\cite[Theorem~1.39]{MR95j:18001}, the subcategory
$\inj_{\rts}(\mathcal{S})$ is then a reflective subcategory of $\rts$.
For any regular transition system $X$, the map $X\to \mathbf{1}$
factors as a composite $X\to F(X) \to \mathbf{1}$ where the left-hand
map belongs to $\cell_{\rts}(\mathcal{S})$ and the right-hand map
belongs to $\inj_{\rts}(\mathcal{S})$ by using the small object
argument in the locally presentable category $\rts$. Then $F(X)$ is
$\mathcal{S}$-orthogonal by Proposition~\ref{simpl}. We deduce that
the map $X\to F(X)$ factors uniquely as a composite $X\to
\bl^{\rts}_{\mathcal{S}}(X)\to F(X)$ by the property of the
adjunction. But the map $X\to \bl^{\rts}_{\mathcal{S}}(X)$ factors
uniquely as a composite $X\to F(X)\to \bl^{\rts}_{\mathcal{S}}(X)$
since the map $X\to F(X)$ belongs to $\cell_{\rts}(\mathcal{S})$ and
since $\bl^{\rts}_{\mathcal{S}}(X)$ is $\mathcal{S}$-orthogonal. Hence
the functor $F$ and $\bl^{\rts}_{\mathcal{S}}$ are isomorphic. \epf

The next proposition compares the functor
$\bl^{\rts}_{\mathcal{S}}:\rts\to\rts$ with the functor
$\bl^{\cts}_{\mathcal{S}}:\cts\to\cts$ defined in an analogous way in
\cite{cubicalhdts}:

\bp \label{comp-L} \label{bl-diff} Let $X$ be a regular transition
system. Then one has the natural isomorphism
\[\CSA_2(\bl^{\cts}_{\mathcal{S}}(X)) \iso \bl^{\rts}_{\mathcal{S}}(X).\]
\ep 

\bpf The map $\bl^{\cts}_{\mathcal{S}}(X)\to
\CSA_2(\bl^{\cts}_{\mathcal{S}}(X))$ is bijective on actions by
Proposition~\ref{calcul_CSA2}. Hence the labelling map of
$\CSA_2(\bl^{\cts}_{\mathcal{S}}(X))$ is one-to-one since the
labelling map of $\bl^{\cts}_{\mathcal{S}}(X)$ is one-to-one by
Proposition~\ref{simpl}. Since the map $X\to \mathbf{1}$ factors as a
composite
\[X\longrightarrow \bl^{\cts}_{\mathcal{S}}(X) \longrightarrow
\CSA_2(\bl^{\cts}_{\mathcal{S}}(X)) \longrightarrow \mathbf{1},\] and
since $\CSA_2(\bl^{\cts}_{\mathcal{S}}(X))$ is $\mathcal{S}$-injective
and regular, the latter satisfies the same universal property as
$\bl^{\rts}_{\mathcal{S}}(X)$. Hence the proof is complete.  \epf

\bth (Compare with \cite[Theorem~8.10]{cubicalhdts}) \label{carw-rts}
A map of regular transition systems $f:X\to Y$ is a weak equivalence
of the Bousfield localization $\bl_{\cub}\rts$ of $\rts$ by the
set of maps $\mathcal{S}$ if and only if the map
$\bl^{\rts}_{\mathcal{S}}(f):\bl^{\rts}_{\mathcal{S}}(X) \to
\bl^{\rts}_{\mathcal{S}}(Y)$ is an isomorphism. \eth

\bpf We already saw in the proof of Theorem~\ref{eg-localizer} that
every map of $\cell_{\rts}(\mathcal{S})$ is a weak equivalence of
$\bl_{\cub}\rts$. This implies that for all morphisms of
regular transition systems $f:X\to Y$, if
$\bl^{\rts}_{\mathcal{S}}(f)$ is an isomorphism, then $f$ belongs to
$\W(\mathcal{S})$. Conversely, let us suppose that $f:X\to Y$ is a
weak equivalence of $\bl_{\cub}\rts$.  Then
$\bl^{\rts}_{\mathcal{S}}(f):\bl^{\rts}_{\mathcal{S}}(X)\to
\bl^{\rts}_{\mathcal{S}}(Y)$ is a map of regular transition systems
between two $\mathcal{S}$-injective regular transition systems. By
Proposition~\ref{simpl}, both $\bl^{\rts}_{\mathcal{S}}(X)$ and
$\bl^{\rts}_{\mathcal{S}}(Y)$ satisfy CSA1 and are
$\mathcal{S}^{cof}$-orthogonal. By
\cite[Proposition~7.7]{cubicalhdts}, both
$\bl^{\rts}_{\mathcal{S}}(X)$ and $\bl^{\rts}_{\mathcal{S}}(Y)$ are
fibrant in $\bl_{\cub}\cts$, and therefore fibrant in
$\bl_{\cub}\rts$ by Proposition~\ref{fibrant-rts}. In other
terms, $\bl^{\rts}_{\mathcal{S}}(f):\bl^{\rts}_{\mathcal{S}}(X)\to
\bl^{\rts}_{\mathcal{S}}(Y)$ is a weak equivalence between two
cofibrant-fibrant objects of the Bousfield localization. Hence,
$\bl^{\rts}_{\mathcal{S}}(f)$ is a weak equivalence of the left
determined model structure $\rts$. By Proposition~\ref{simpl}, both
$\bl^{\rts}_{\mathcal{S}}(X)$ and $\bl^{\rts}_{\mathcal{S}}(Y)$
satisfy CSA1.  By Proposition~\ref{homotopy_is_equality}, one deduces
that $\bl^{\rts}_{\mathcal{S}}(f)$ is an isomorphism.  \epf

Proposition~\ref{compare-cellS} and Theorem~\ref{th-compare-weak} help
to understand the difference between the weak equivalences of
$\bl_{\cub}\cts$ and of $\bl_{\cub}\rts$.

\bp \label{compare-cellS} For all cubical transition systems $X$, the
map $X\to \bl^{\cts}_{\mathcal{S}}(X)$ is bijective on states and onto
on actions. There exists a cubical transition system $X_0$ such that
the map $X_0\to \bl^{\cts}_{\mathcal{S}}(X_0)$ is not onto on
transitions.  For all regular transition systems $Y$, the map $Y\to
\bl^{\rts}_{\mathcal{S}}(Y)$ is onto on states, on actions and on
transitions.  There exists a regular transition system $Y_0$ such that
the map $Y_0\to \bl^{\rts}_{\mathcal{S}}(Y_0)$ is not bijective on
states. \ep

\bpf This is a corollary of Proposition~\ref{surjectivity-cellS-cts}
and Proposition~\ref{surjectivity-cellS} of Appendix~\ref{cellS}. \epf

\bth \label{th-compare-weak} There exists a strict inclusion of sets
\begin{multline*}\left\{\hbox{weak equivalences of $\bl_{\cub}\cts$ between regular
      t.s.}\right\} \\\subset \left\{\hbox{weak equivalences of
      $\bl_{\cub}\rts$}\right\}.\end{multline*} In other terms, if
$f:X\to Y$ is a weak equivalence of $\bl_{\cub}\cts$ between two
regular transition systems, then $f$ is a weak equivalence of
$\bl_{\cub}\rts$.  There exists a weak equivalence of $\bl_{\cub}\rts$
which is not a weak equivalence of $\bl_{\cub}\cts$.  \eth

\bpf Let $f:X\to Y$ be a weak equivalence of $\bl_{\cub}\cts$ between
two regular transition systems. Then by
\cite[Theorem~8.10]{cubicalhdts}, the map
$\bl^{\cts}_{\mathcal{S}}(f)$ is an isomorphism. The map
$\CSA_2(\bl^{\cts}_{\mathcal{S}}(f))$ is therefore an isomorphism. So,
by Proposition~\ref{comp-L}, $\bl^{\rts}_{\mathcal{S}}(f)$ is an
isomorphism.  Hence by Theorem~\ref{carw-rts}, $f$ is a weak
equivalence of $\bl_{\cub}\rts$.

Now we want to find a weak equivalence $g$ of $\bl_{\cub}\rts$ which
is not a weak equivalence of $\bl_{\cub}\cts$.  One has $\omega
(C_2[x,x]) = (\{0,1\}^2,\{(x,1),(x,2)\})$ by
Proposition~\ref{cas_cube} with $x\in \Sigma$.  Consider the set
$\{0,1\}^2\p\{-,+\}$ and let us make the identifications
$(0,0,-)=(0,0,+)=I$ and $(1,1,-)=(1,1,+)=F$. Write $S$ for the
quotient. Let $W=(S,\{u,v^-,v^+\})$. For $\alpha\in \{-,+\}$, consider
the map $\phi^\alpha:\omega (C_2[x,x])\to W$ of $\set^{\{s\}\cup
  \Sigma}$ induced by the mappings $(\epsilon_1,\epsilon_2)\mapsto
(\epsilon_1,\epsilon_2,\alpha)$ for $(\epsilon_1,\epsilon_2)\in
\{0,1\}^2$, $(x,1)\mapsto u$ and $(x,2)\mapsto v^\alpha $.  Consider
the $\omega$-final lift $\overline{W}$ of the cone of maps
$\phi^-,\phi^+:\omega (C_2[x,x]) \rightrightarrows W$.  By
Theorem~\ref{cubical-lift}, the weak transition system $\overline{W}$
is cubical. The only higher dimensional transitions of $\overline{W}$
are the four transitions $(I,u,v^\pm,F)$ and $(I,v^\pm,u,F)$. Hence
the unique state $\nu$ such that the tuples $(I,u,\nu)$ and
$(\nu,v^\pm,F)$ are transitions of $\overline{W}$ is
$\nu=(1,0,\pm)$. It turns out that the unique state $\nu'$ such that
the tuples $(I,v^\pm,\nu')$ and $(\nu',u,F)$ are transitions of
$\overline{W}$ is $\nu'=(0,1,\pm)$. One deduces that $\overline{W}$ is
regular. There exists a map of cubical transition systems
$g:\overline{W} \to C_2[x,x]$ defined as follows: it takes the state
$(\epsilon_1,\epsilon_2,\pm)$ to $(\epsilon_1,\epsilon_2)$ for
$(\epsilon_1,\epsilon_2)\in, \{0,1\}^2$, the action $u$ to $(x,1)$ and
the actions $v^-$ and $v^+$ to $(x,2)$. It is easy to see that one has
the isomorphisms
\[\bl^{\rts}_{\mathcal{S}}(\overline{W}) \iso C_2[x,x] \iso \bl^{\rts}_{\mathcal{S}}(C_2[x,x]),\]
hence $g$ is a weak equivalence of $\bl_{\cub}\rts$ by
Theorem~\ref{carw-rts}. Since $g$ is not bijective on states, the map
$\bl^{\cts}_{\mathcal{S}}(f)$ is not bijective on states by
Proposition~\ref{compare-cellS}. Therefore the map
$\bl^{\cts}_{\mathcal{S}}(f)$ is not an isomorphism. Hence $g$ is not
a weak equivalence of $\bl_{\cub}\cts$ by
\cite[Theorem~8.10]{cubicalhdts}.  \epf

We can now completely elucidate this model structure thanks to the
following result:

\bth (Compare with \cite[Theorem~8.11 (4)(5)]{cubicalhdts}) The left
adjoint \[\bl^{\rts}_{\mathcal{S}}:\bl_{\cub}\rts \to
\inj_{\rts}(\mathcal{S})\] induces a left Quillen equivalence between
\[\bl_{\cub}\rts \to \inj_{\rts}(\mathcal{S})\] equipped with
the discrete model structure (all maps are cofibrations and fibrations
and the weak equivalences are the isomorphisms).  \eth

\bpf For any fibrant object $X$ of $\inj_{\rts}(\mathcal{S})$, the map
$\bl^{\rts}_{\mathcal{S}}(X) \to X$ is an isomorphism and $X$ is
cofibrant in $\bl_{\cub}\rts$. For any cofibrant object $Y$ of
$\bl_{\cub}\rts$, $Y$ is fibrant in $\inj_{\rts}(\mathcal{S})$
and the map $Y\to \bl^{\rts}_{\mathcal{S}}(Y)$ is a weak equivalence
of $\bl_{\cub}\rts$ by Proposition~\ref{def-LS-rts} and by
Lemma~\ref{for-when-cell-is-weak-2}. This is the definition of a
Quillen equivalence. \epf

Theorem~\ref{carw-rts} does not mean that two regular transition
systems are weakly equivalent if and only if they are
isomorphic. Indeed, for any regular transition system $X$, the unit
map $X\to \bl^{\rts}_{\mathcal{S}}(X)$, by identifying the actions of
$X$ with their labelling, modifies the geometric structure of $X$ by
forcing identifications of states (see
Proposition~\ref{compare-cellS}). Roughly speaking, this map removes
all non-discernable transitions. This behaviour is slightly different
from the one of the unit map $X\to \bl^{\cts}_{\mathcal{S}}(X)$. Once
again by Proposition~\ref{compare-cellS}, the unit map $X\to
\bl^{\cts}_{\mathcal{S}}(X)$ also identifies the actions of a cubical
transition system $X$ by their labelling, but the latter map is
constant on states, and not necessarily onto on transitions. It may
create new transitions which are actually not observable and which are
killed by applying the functor $\CSA_2:\cts\to \rts$.

\section{Fibrant regular and cubical transition systems}
\label{fibrant-char}

The purpose of this last section is to describe completely the fibrant
regular and cubical transition systems.  We already know by
Proposition~\ref{fibrant-rts} that the fibrant regular transition
systems are exactly the fibrant cubical ones which are regular. Thus,
we just have to give a combinatorial characterization of the fibrant
objects of $\bl_{\cub}\cts$.  Corollary~\ref{car-fibrant-lscts}
encompasses the results of \cite{cubicalhdts} and
\cite{erratum_cubicalhdts}.

\bd A cubical transition system $X$ is {\rm combinatorially fibrant}
if for any $n\geq 1$, any state $\alpha$ and $\beta$ and any actions
$u_1,v_1,\dots,u_n,v_n$ such that $\mu(u_i)=\mu(v_i)$ for $1\leq i
\leq n$, if the tuple $(\alpha,u_1,\dots,u_n,\beta)$ is a transition
of $X$, then the tuple $(\alpha,v_1,\dots,v_n,\beta)$ is a transition
of $X$ as well. \ed

\bp \label{calculation-path} Let $X=(S,\mu:L\to \Sigma,T)$ be a
combinatorially fibrant cubical transition system. Write
$\cocyl:\cts\to\cts$ for the right adjoint of the cartesian cylinder
$\cyl:\cts\to\cts$.  Then the cubical transition system $\cocyl(X)$
has $S$ as its set of states and $L\p_\Sigma L$ as its set of actions,
the labelling map is the composite map $\mu:L\p_\Sigma L \to L \to
\Sigma$ and a tuple $(\alpha,(u^0_1,u^1_1),\dots,(u^0_n,u^1_n),\beta)$
of $S\p (L\p_\Sigma L)^n \p S$ is a transition of $\cocyl(X)$ if and
only if there exist $\epsilon_1,\dots,\epsilon_n \in \{0,1\}$ such
that the tuple
$(\alpha,u^{\epsilon_1}_1,\dots,u^{\epsilon_n}_n,\beta)$ is a
transition of $X$. \ep

\bpf Let us recall that the cartesian cylinder $\cyl:\cts \to \cts$ is
the restriction of an endofunctor of $\wts$ defined in the same
way. The functor $\cyl:\wts\to \wts$ has a right adjoint
$\cocyl^{\wts}:\wts \to \wts$ defined on objects as follows
\cite[Proposition~5.8]{cubicalhdts}: for a weak transition system
$X=(S,\mu:L\to \Sigma,T)$, the weak transition system
$\cocyl^{\wts}(X)$ has the same set of states $S$, the set of actions
is $L\p_\Sigma L$ and a tuple $(\alpha,(u^-_1,u^+_1), \dots,
(u^-_n,u^+_n), \beta)$ with $n\geq 1$ is a transition of
$\cocyl^{\wts}(X)$ if and only if the $2^n$ tuples $(\alpha, u_1^\pm,
\dots, u_n^\pm, \beta)$ are transitions of $X$. The right adjoint of
the functor $\cyl:\cts \to \cts$ is equal to the composite
functor \[\xymatrix{ \cocyl : \cts \subset \wts \fR{\cocyl^{\wts}} &&
  \wts \fR{} && \cts,}\] where the right-hand functor from $\wts$ to
$\cts$ is the coreflection.

Let $(u,v)\in L\p_\Sigma L$. Since $u$ is used in $X$, there exists a
transition $(\alpha,u,\beta)$ of $X$. Since $\mu(u)=\mu(v)$ and since
$X$ is combinatorially fibrant, the triple $(\alpha,v,\beta)$ is a
transition of $X$. This means that the couple $(u,v)\in L\p_\Sigma L$
is used by the transition $(\alpha,(u,v),\beta)$ of
$\cocyl^{\wts}(X)$. We deduce that all actions of $\cocyl^{\wts}(X)$
are used. Consider a transition \[(\alpha,(u^-_1,u^+_1), \dots,
(u^-_n,u^+_n), \beta)\] of $\cocyl^{\wts}(X)$ with $n\geq 2$. Let
$1\leq p \leq n-1$. Since $X$ is cubical, there exists a state
$\gamma$ such that the tuples $(\alpha, u^-_1,\dots,u^-_p,\gamma)$ and
$(\gamma,u^-_{p+1},\dots,u^-_n,\beta)$ are two transitions of $X$. But
$X$ is combinatorially fibrant. This implies that all tuples $(\alpha,
u^\pm_1,\dots,u^\pm_p,\gamma)$ and
$(\gamma,u^\pm_{p+1},\dots,u^\pm_n,\beta)$ are transitions of
$X$. Therefore the two tuples
\[
(\alpha,(u^-_1,u^+_1), \dots, (u^-_p,u^+_p),
\gamma),(\gamma,(u^-_{p+1},u^+_{p+1}), \dots, (u^-_n,u^+_n), \beta)
\] are transitions of $\cocyl^{\wts}(X)$. This means that the weak
transition system $\cocyl^{\wts}(X)$ satisfies the Intermediate state
axiom.  We have just proved that if $X$ is combinatorially fibrant,
then the weak transition system $\cocyl^{\wts}(X)$ is cubical: in
other terms, one has $\cocyl(X)=\cocyl^{\wts}(X)$ in this case.
Finally and because $X$ is combinatorially fibrant, all tuples
$(\alpha, u_1^\pm, \dots, u_n^\pm, \beta)$ are transitions of $X$ if
and only if there exist $\epsilon_1,\dots,\epsilon_n \in \{0,1\}$ such
that the tuple
$(\alpha,u^{\epsilon_1}_1,\dots,u^{\epsilon_n}_n,\beta)$ is a
transition of $X$. This completes the proof.  \epf

\bp \label{fact0} If the cubical transition system $X$ is
combinatorially fibrant, then so is the cubical transition system
$\cocyl(X)$. \ep

\bpf Let $X=(S,\mu:L\to \Sigma,T)$ be a combinatorially fibrant
cubical transition system. Let \[(\alpha,(u^-_1,u^+_1), \dots,
(u^-_n,u^+_n), \beta),(\alpha,(v^-_1,v^+_1), \dots, (v^-_n,v^+_n),
\beta)\] be two tuples of $S\p (L\p_\Sigma L)^n \p S$ with $n\geq 1$
and $\mu(u^-_i,u^+_i)=\mu(v^-_i,v^+_i)$ for $1\leq i \leq n$. Let us
suppose that $(\alpha,(u^-_1,u^+_1), \dots, (u^-_n,u^+_n), \beta)$ is
a transition of $\cocyl(X)$. Then the tuple
$(\alpha,u^-_1,\dots,u^-_n,\beta)$ is a transition of $X$.  But for
all $1\leq i \leq n$, one has
\[\mu(u^-_i)=\mu(u^+_i)=\mu(u^-_i,u^+_i)=\mu(v^-_i,v^+_i)=\mu(v^-_i)=\mu(v^+_i).\]
So, all tuples $(\alpha, v_1^\pm, \dots, v_n^\pm, \beta)$ are
transitions of $X$ because $X$ is combinatorially fibrant. This
implies that the tuple $(\alpha,(v^-_1,v^+_1), \dots, (v^-_n,v^+_n),
\beta)$ is a transition of $\cocyl(X)$. This is the definition of
combinatorial fibrancy applied to $\cocyl(X)$.  \epf

\bp \label{fact1} Let $X$ be a cubical transition system. If $X$ is
combinatorially fibrant, then it is injective with respect to any map
of the form $f\star \gamma^\epsilon$ for $\epsilon=0,1$ for any
cofibration of cubical transition systems $f$. \ep

\bpf Let $f:A\to B$ be a map of cubical transition systems. Let $L$ be
the set of actions of $X$. By adjunction, the cubical transition
system $X$ is injective with respect to $f\star \gamma^\epsilon$ if
and only if the map $\pi^\epsilon:\cocyl(X)\to X$ satisfies the RLP
with respect to $f$. Let us recall that the map
$\pi^\epsilon:\cocyl(X)\to X$ is the identity on states and the
projection on the $(\epsilon+1)$-th component $L\p_\Sigma L \to L$ on
actions by Proposition~\ref{calculation-path}. Consider a diagram of
solid arrows of cubical transition systems:
\[
\xymatrix
{
A \fD{f} \fR{\phi} && \cocyl(X) \fD{\pi^\epsilon} \\
&& \\
B \fR{\psi}\ar@{-->}[rruu]^-{\ell} && X.
}
\]
Since the right vertical map is onto on actions and the left vertical
map is one-to-one on actions, there exists a set map
$\widetilde{\ell}:L_B \to L\p_\Sigma L$ from the set of actions of
$B$ to the set of actions of $\cocyl(X)$ such that the following
diagram of sets is commutative, $L_A$ being the set of actions of $A$
(note that $\widetilde{\pi^\epsilon}$ is the projection on the
$(\epsilon+1)$-th component):
\[
\xymatrix
{
L_A \fD{f} \fR{\phi} && L\p_\Sigma L \fD{\pi^\epsilon} \\
&& \\
L_B \fR{\psi}\ar@{-->}[rruu]^-{\widetilde{\ell}} && L.
}
\]
Let $\ell:B\to \cocyl(X)$ defined on states by
$\ell(\alpha)=\psi(\alpha)$ and on actions by
$\ell(u)=\widetilde{\ell}(u)$. The diagram
\[
\xymatrix
{
A \fD{f} \fR{\phi} && \cocyl(X) \fD{\pi^\epsilon} \\
&& \\
B \fR{\psi}\ar@{->}[rruu]^-{\ell} && X
}
\]
is commutative since its right vertical map is the identity on
states. It just remains to prove that $\ell:B\to \cocyl(X)$ is a
well-defined map of cubical transition systems. Let
$(\alpha,u_1,\dots,u_n,\beta)$ be a transition of $B$. It suffices to
prove that the tuple
\[(\alpha,\widetilde{\ell}(u_1),\dots,\widetilde{\ell}(u_n),\beta)\] is
a transition of $\cocyl(X)$ to complete the proof. Without lack of
generality, we can suppose that $\epsilon=0$, which means that
$\widetilde{\ell}(u)=(\psi(u),\chi(u))$. One obtains
\[ (\alpha,\widetilde{\ell}(u_1),\dots,\widetilde{\ell}(u_n),\beta) =
(\alpha,(\psi(u_1),\chi(u_1)),\dots,(\psi(u_n),\chi(u_n)),\beta).\]
Since $\psi$ maps the transitions of $B$ to transitions of $X$, the
tuple $(\alpha,\psi(u_1),\dots,\psi(u_n),\beta)$ is a transition of
$X$.  Since $\mu(\psi(u))=\mu(u)=\mu(\chi(u))$ for all actions $u$ of
$B$, and since $X$ is combinatorially fibrant, the tuple
$(\alpha,(\psi(u_1),\chi(u_1)),\dots,(\psi(u_n),\chi(u_n)),\beta)$ is
then a transition of $\cocyl(X)$ by
Proposition~\ref{calculation-path}.  \epf

\bp\label{fact2} Let $X$ be a cubical transition system. If $X$ is
combinatorially fibrant, then it is injective with respect to the maps
of $\mathcal{S}^{cof}$.  \ep

\bpf Let $x\in \Sigma$. Consider a diagram of solid arrows of cubical
transition systems 
\[
\xymatrix
{
C_1[x]\sqcup C_1[x] \fD{p_x^{cof}} \fR{\phi} && X \\
&& \\
Z^{x_1,x_2}_x, \ar@{-->}[rruu]^-{\ell}&&
}
\]
where $x_1$ and $x_2$ are the two actions of $C_1[x]\sqcup C_1[x]$ and
where $Z^{x_1,x_2}_x$ is the cubical transition system depicted in
Figure~\ref{Z^{x_1,x_2}_x}.  Define $\ell$ on states by
$\ell(\alpha)=\phi(\alpha)$, and on actions by $\ell(x_i)=\phi(x_i)$
for $i=1,2$ and $\ell(x)=\phi(x_1)$.  Let
$(\alpha_i,\phi(x_i),\beta_i)$ for $i=1,2$ be the images by $\phi$ of
the two transitions of $C_1[x]\sqcup C_1[x]$. Since $X$ is
combinatorially fibrant, the two triples
$(\alpha_i,\phi(x_{3-i}),\beta_i)$ for $i=1,2$ are two transitions of
$X$. The map $\ell$ is therefore a well-defined map of cubical
transition systems.  \epf

\bp \label{produit_cts} Let $X = (S,\mu:L \rightarrow \Sigma, T)$ and
$X' = (S',\mu':L' \rightarrow \Sigma, T')$ be two cubical transition
systems. The binary product $X\p X'$ has $S\p S'$ as its set of
states, $L\p_\Sigma L' = \{(x,x')\in L\p L', \mu(x) = \mu'(x')\}$ as
its set of actions and the labelling map $\mu\p_\Sigma \mu':L\p_\Sigma
L' \rightarrow \Sigma$. A tuple $((\alpha,\alpha'), (u_1,u'_1), \dots,
(u_n,u'_n), (\beta,\beta'))$ is a transition of $X\p X'$ if and only
if $\mu(u_i) = \mu'(u'_i)$ for $1\leq i \leq n$ with $n\geq 1$, the
tuple $(\alpha, u_1, \dots, u_n, \beta)$ is a transition of $X$ and
$(\alpha', u'_1, \dots, u'_n, \beta')$ a transition of $X'$. \ep

\bpf The binary product is the same in $\cts$ and in $\wts$ because
$\cts$ is a small-injectivity class of $\wts$. The theorem is then a
consequence of \cite[Proposition~5.5]{cubicalhdts}.  \epf

\bp \label{fact3} Let $X$ be a cubical transition system. If $X$ is
combinatorially fibrant, then it is injective with respect to any map
of the form $f\star \gamma$ for any map of cubical transition systems
$f$ which is onto on states. \ep

\bpf Let $f:A\to B$ be a map of cubical transition systems. By
adjunction, the cubical transition system $X$ is injective with
respect to $f\star \gamma$ if and only if the map $\pi:\cocyl(X)\to X
\p X$ satisfies the RLP with respect to $f$.  Consider a diagram of
solid arrows of cubical transition systems:
\[
\xymatrix
{
A \fD{f} \fR{\phi} && \cocyl(X) \fD{\pi} \\
&& \\
B \fR{\psi=(\psi_0,\psi_1)}\ar@{-->}[rruu]^-{\ell} && X\p X.
}
\]
Since the set map $f:A^0\to B^0$ is onto by hypothesis, for any state
$\alpha$ of $B$, there exists $s(\alpha)\in A^0$ such that
$f(s(\alpha))=\alpha$.  Let $\ell:B\to \cocyl(X)$ defined on states by
$\ell(\alpha)=\phi(s(\alpha))$ and on actions by $\ell(u)=\psi(u)$
(since $X$ is combinatorially fibrant, the map $\pi:\cocyl(X)\to X\p
X$ is the identity on actions by Proposition~\ref{calculation-path}).
We are going to prove that $\ell$ is a well-defined map of cubical
transition systems and that it is a lift of the diagram above.

\underline{$\ell$ is a lift for the sets of actions.}  One has the
following diagram of solid arrows between the sets of actions:
\[
\xymatrix
{
L_A \fD{f} \fR{\phi} && L_X \p_\Sigma L_X \ar@{=}[dd] \\
&& \\
L_B \fR{\psi}\ar@{-->}[rruu]^-{\psi}  && L_X \p_\Sigma L_X.
}
\]
It is evident that the two triangles commute since the square of solid
arrows commutes.

\underline{$\ell$ is a lift for the sets of states.}  One has the
diagram of solid arrows between the sets of states:
\[
\xymatrix
{
A^0 \fD{f} \fR{\phi} && X^0 \ar@{->}[dd]^-{\Delta} \\
&& \\
B^0 \fR{\psi=(\psi_0,\psi_1)}\ar@{-->}[rruu]^-{\phi\circ s} \ar@{-->}@/^20pt/[uu]^-{s} && X^0 \p X^0,
}
\]
where $\Delta:s\mapsto (s,s)$ is the codiagonal map.   For any state
$\beta$ of $B^0$, one has
\begin{align*}
\psi(\beta) &= \psi(f(s(\beta)))& \hbox{ since $s$ is a section of $f$}\\
&=\pi(\phi(s(\beta))) & \hbox{since $\psi\circ f=\pi\circ \phi$}\\
&=(\phi(s(\beta)),\phi(s(\beta))) & \hbox{by Proposition~\ref{produit_cts}.}
\end{align*}
Hence we obtain $\psi_0=\psi_1=\phi\circ s$ on states, and therefore
$\Delta\circ \phi\circ s=\psi$ on states. We deduce that the bottom
triangle commutes on states.  For any state $\alpha$ of $A^0$, one has
\begin{align*}
\Delta(\phi(s(f(\alpha)))) &=\psi(f(s(f(\alpha))))& \hbox{ since $\Delta\circ\phi=\psi\circ f$}\\
&=\psi(f(\alpha)) & \hbox{ since $s$ is a section of $f$}\\
&=\Delta(\phi(\alpha))& \hbox{because the square above is commutative.}
\end{align*}
Hence we obtain $\phi\circ s\circ f=\phi$ on states. We obtain that
the top triangle commutes.

\underline{$\ell$ maps a transition of $B$ to a transition of
  $\cocyl(X)$.}  Let $(\alpha,u_1,\dots,u_n,\beta)$ be a transition of
$B$. Then one has
\begin{multline*}
(\ell(\alpha),\ell(u_1),\dots,\ell(u_n),\ell(\beta)) =
(\phi(s(\alpha)),\psi(u_1),\dots,\psi(u_n),\phi(s(\beta))) \\ = 
(\psi_0(\alpha),\psi(u_1),\dots,\psi(u_n),\psi_0(\beta)).
\end{multline*}
The tuple
$(\psi_0(\alpha),\psi_0(u_1),\dots,\psi_0(u_n),\psi_0(\beta))$ is a
transition of $X$ since it is the image by the composite map of
cubical transition systems $\psi_0:B\to X\p X \to X$ of the transition
$(\alpha,u_1,\dots,u_n,\beta)$ of $B$. Therefore by
Proposition~\ref{calculation-path} applied with
$\epsilon_1=\dots=\epsilon_n=0$, the tuple
\[(\ell(\alpha),\ell(u_1),\dots,\ell(u_n),\ell(\beta))\] is a
transition of $\cocyl(X)$ since $X$ is combinatorially fibrant. This
means that \[\ell:B\longrightarrow \cocyl(X)\] is a well-defined map
of cubical transition systems. \epf

\bp \label{fact4} Let $X$ be a cubical transition system. If $X$ is
combinatorially fibrant, then it is injective with respect to any map
of the form $(f\star \gamma) \star \gamma $ for any map of cubical
transition systems $f$.  \ep

\bpf Let $f:A\to B$ be a map of cubical transition systems. The map
$f\star \gamma$ goes from $(B\sqcup B) \sqcup_{A\sqcup A} \cyl(A)$ to
$\cyl(B)$.  Since the forgetful functor from $\cts$ to $\set$ taking a
cubical transition system to its underlying set of states is
colimit-preserving, the set of states of the source of $f\star \gamma$
is $B^0 \sqcup_{A^0} B^0$.  Hence the map $f\star \gamma$ is onto on
states. Then by Proposition~\ref{fact3}, $X$ is injective with respect
to $(f\star \gamma) \star \gamma$. \epf

\begin{nota} Let $I$ and $S$ be two sets of maps of a locally
  presentable category $\K$. Let $\cyl:\K\to \K$ be a cylinder. Denote
  by $\Lambda_\K(\cyl,S,I)$ the set of maps defined as follows:
\begin{itemize}
\item $\Lambda^0_\K(\cyl,S,I) = S \cup (I \star \gamma^0) \cup (I
  \star \gamma^1)$
\item $\Lambda^{n+1}_\K(\cyl,S,I) =
  \Lambda^n_\K(\cyl,S,I) \star \gamma$
\item $\Lambda_\K(\cyl,S,I) = \bigcup_{n\geq 0}
  \Lambda^n_\K(\cyl,S,I)$.
\end{itemize} 
\end{nota}

\bth \label{sens1} Let $X$ be a cubical transition system. If $X$ is
combinatorially fibrant, then it is fibrant. \eth

\bpf Let $X$ be a combinatorially fibrant cubical transition
system. By Proposition~\ref{fact1} and Proposition~\ref{fact2}, it is
$\Lambda^0(\cyl,\mathcal{S}^{cof},\mathcal{I})$-injective.  Let
$f:A\to B$ be a map of cubical transition systems. Let $\epsilon\in
\{0,1\}$.  The map $f\star \gamma^\epsilon$ goes from $B\sqcup_A
\cyl(A)$ to $\cyl(B)$.  Since the forgetful functor from $\cts$ to
$\set$ taking a cubical transition system to its underlying set of
states is colimit-preserving, the set of states of the source of
$f\star \gamma^\epsilon$ is $B^0$.  Hence $f\star \gamma^\epsilon$ is
bijective on states. Therefore all maps of
$\Lambda^0(\cyl,\mathcal{S}^{cof},\mathcal{I})$ are bijective on
states. Then, by Proposition~\ref{fact3}, $X$ is
$\Lambda^1(\cyl,\mathcal{S}^{cof},\mathcal{I})$-injective.  The
cubical transition system $X$ is
$\Lambda^n(\cyl,\mathcal{S}^{cof},\mathcal{I})$-injective for all
$n\geq 2$ by Proposition~\ref{fact4}. Hence $X$ is fibrant in the
Bousfield localization of $\cts$ by the cofibrations of
$\mathcal{S}^{cof}$ by \cite[Corollary~6.8]{cubicalhdts} and
\cite[Theorem~4.6]{homotopyprecubical}. But Bousfield localizing by
$\mathcal{S}^{cof}$ is the same as Bousfield localizing by
$\mathcal{S}$, which is the same as Bousfield localizing by the
cubification functor. Hence the proof is complete. \epf

\begin{nota} Let $x\in \Sigma$. The two maps from $C_1[x]$ to $\dd{x}$
  are denoted by $c_x^\epsilon$ for $\epsilon=0,1$. One has $p_x=c_x^0
  \sqcup c_x^1$ for all $x\in \Sigma$. \end{nota}

\bp \label{lc1} Let $x\in \Sigma$.  Consider the pushout
  diagram of $\cts$
\[
\xymatrix
{
C_1[x] \fR{c_x^0} \fD{\gamma^0_{C_1[x]}} && \dd{x} \fD{} \\ 
&& \\
\cyl(C_1[x])  \fR{} && \cyl(C_1[x]) \sqcup_{0,0} \dd{x}. \cocartesien
}
\]
The composite
\[
\theta_x:\xymatrix{C_1[x] \sqcup C_1[x] \fR{\gamma^1_{C_1[x]}\sqcup
  c_x^1} && \cyl(C_1[x]) \sqcup  \dd{x} \fR{} && \cyl(C_1[x]) \sqcup_{0,0} \dd{x}}
\] 
is a trivial cofibration of $\bl_{\cub}\cts$.
\ep

\bpf The map $\theta_x$ is depicted in Figure~\ref{theta}.  It is
bijective on actions, therefore it is a cofibration. One has
$\bl^{\cts}_{\mathcal{S}}(C_1[x] \sqcup C_1[x]) \iso
\bl^{\cts}_{\mathcal{S}}(\cyl(C_1[x]) \sqcup_{0,0} \dd{x}) \iso
\dd{x}$. Hence it is a weak equivalence of $\bl_{\cub}\cts$ by
\cite[Theorem~8.10]{cubicalhdts}.  \epf

\begin{figure}
\[\left\{\begin{array}{c} C_1[x] \sqcup C_1[x] \\
        {
\xymatrix{
\alpha \fR{x_1} && \beta \\
\gamma  \fR{x_2} && \delta
}
}\end{array}\right.
  \stackrel{\theta_{x}}\longrightarrow \left\{\begin{array}{c} \cyl(C_1[x]) \sqcup_{0,0} \dd{x} \\ 
{\xymatrix{
\alpha \ar@/^5pt/[rr]^-{x_1} \ar@/_5pt/[rr]_-{x_2}&& \beta \\ \gamma  \fR{x_2} && \delta\\ }}\end{array}\right.\]
\caption{Cofibration $\theta_x$ with $\mu(x_1) =
  \mu(x_2) = x$}
\label{theta}
\end{figure}

\bp \label{lc2} In the following, the notation
$\sqcup_{\tiny\begin{array}{c}0_n=0_n\\1_n=1_n\end{array}}$ means the
identification of the initial states (the final states resp.) of the
two summands. Let $n\geq 2$ and $x_1,\dots,x_n \in \Sigma$. Then the
map
\begin{multline*}\eta_{x_1,\dots,x_n}:\de C_n[x_1,\dots,x_n] \sqcup_{\tiny\begin{array}{c}0_n=0_n\\1_n=1_n\end{array}} C_n[x_1,\dots,x_n] 
  \\ \longrightarrow C_n[x_1,\dots,x_n]
  \sqcup_{\tiny\begin{array}{c}0_n=0_n\\1_n=1_n\end{array}}
  C_n[x_1,\dots,x_n]\end{multline*} induced by the inclusion \[\de
C_n[x_1,\dots,x_n] \subset C_n[x_1,\dots,x_n]\] is a trivial
cofibration of $\bl_{\cub}\cts$.  \ep

\bpf The map $\eta_{x_1,\dots,x_n}$ is bijective on actions: the set
of actions is $\{(x_1,1),\dots,(x_n,n)\}\p \{0,1\}$, with for example
$0$ for the left-hand term and $1$ for the right-hand term.  Hence it
is a cofibration. The map $\eta_{x_1,\dots,x_n}$ is also bijective on
states: the set of states is a set denoted by $\{0,1\}^n
\sqcup_{\tiny\begin{array}{c}0_n=0_n\\1_n=1_n\end{array}} \{0,1\}^n$,
which means the quotient of the coproduct $\{0,1\}^n \sqcup \{0,1\}^n$
by the identifications of $0_n$ ($1_n$ resp.) of the left-hand term
with $0_n$ ($1_n$ resp.) of the right-hand term. Since the map $X\to
\bl^{\cts}_{\mathcal{S}}(X)$ is bijective on states for all cubical
transition systems $X$, the map of cubical transition systems
\begin{multline*} 
\bl^{\cts}_{\mathcal{S}}(\eta_{x_1,\dots,x_n}):\bl^{\cts}_{\mathcal{S}}\left(\de C_n[x_1,\dots,x_n] \sqcup_{\tiny\begin{array}{c}0_n=0_n\\1_n=1_n\end{array}} C_n[x_1,\dots,x_n]\right) \longrightarrow \\ \bl^{\cts}_{\mathcal{S}}\left(C_n[x_1,\dots,x_n]
\sqcup_{\tiny\begin{array}{c}0_n=0_n\\1_n=1_n\end{array}}
C_n[x_1,\dots,x_n]\right)
\end{multline*} 
is bijective on states as well. The set of actions of the source and
target of $\bl^{\cts}_{\mathcal{S}}(\eta_{x_1,\dots,x_n})$ is
$\{x_1,\dots,x_n\}$. Since
$\bl^{\cts}_{\mathcal{S}}(\eta_{x_1,\dots,x_n})$ is one-to-one on
action by \cite[Remark~8.8]{cubicalhdts}, it is bijective on actions.
By \cite[Proposition~4.4]{homotopyprecubical}, the map
$\bl^{\cts}_{\mathcal{S}}(\eta_{x_1,\dots,x_n})$ is one-to-one on
transitions by.  To see that the map
$\bl^{\cts}_{\mathcal{S}}(\eta_{x_1,\dots,x_n})$ is also onto on
transitions, it suffice to see that the $n!$ $n$-transitions of the
left-hand $n$-cube of the target are the $n!$ tuples
$(0_n,x_{\sigma(1)},\dots,x_{\sigma(n)},1_n)$ which are actually
transitions of the source because of the identifications of the two
initial states and the two final states.  So
$\bl^{\cts}_{\mathcal{S}}(\eta_{x_1,\dots,x_n})$ is an
isomorphism. Therefore by \cite[Theorem~8.10]{cubicalhdts}, the map
$\eta_{x_1,\dots,x_n}$ is a weak equivalence of $\bl_{\cub}\cts$.
\epf

\bp \label{combinatorially-fibrant-injective} A cubical transition
system is combinatorially fibrant if and only if it is injective with
respect to $\theta_x$ and $\eta_{x_1,\dots,x_n}$ for all
$x,x_1,\dots,x_n\in \Sigma$.  \ep

\bpf Let $X$ a combinatorially fibrant cubical transition system. Then
$X$ is fibrant by Theorem~\ref{sens1}. Since the maps $\theta_x$ and
$\eta_{x_1,\dots,x_n}$ for all $x,x_1,\dots,x_n\in \Sigma$ are trivial
cofibrations by Proposition~\ref{lc1} and Proposition~\ref{lc2}, $X$
is injective with respect to these maps. Conversely, let $X$ be a
cubical transition system which is injective with respect to
$\theta_x$ and $\eta_{x_1,\dots,x_n}$ for all $x,x_1,\dots,x_n\in
\Sigma$. Let $(\alpha,x_1,\beta)$ be a transition of $X$ and let $x_2$
an action of $X$ such that $\mu(x_1)=\mu(x_2)$.  The injectivity of
$X$ with respect to $\theta_{\mu(x_1)}$ proves that the triple
$(\alpha,x_2,\beta)$ is a transition of $X$. Let
$(\alpha,x_1,\dots,x_n,\beta)$ be a transition of $X$ with $n\geq
2$. Let $y_1,\dots,y_n$ be $n$ actions of $X$ with $\mu(x_i)=\mu(y_i)$
for $1\leq i \leq n$. The injectivity of $X$ with respect to
$\eta_{\mu(x_1),\dots,\mu(x_n)}$ proves that the triple
$(\alpha,y_1,\dots,y_n,\beta)$ is a transition of $X$. So, $X$ is
combinatorially fibrant.  \epf

\begin{cor} \label{sens2} Let $X$ be a cubical transition system. If
  $X$ is fibrant, then it is combinatorially fibrant.  \end{cor}

\bpf Let $X$ be a fibrant cubical transition system. Then it is
injective with respect to any trivial cofibration of
$\bl_{\cub}\cts$. By Proposition~\ref{lc1}, Proposition~\ref{lc2} and
Proposition~\ref{combinatorially-fibrant-injective}, it is then
combinatorially fibrant.  \epf

\begin{cor} \label{car-fibrant-lscts} A cubical transition system $X$
  is fibrant in $\bl_{\cub}\cts$ if and only it is
  combinatorially fibrant. \end{cor}

\begin{cor} \label{Sinj-fibrant}
  Every $\mathcal{S}$-injective cubical transition system is fibrant
  in $\bl_{\cub}\cts$.
\end{cor}

\bpf Let $(\alpha,u_1,\dots,u_n,\beta)$ and
$(\alpha,v_1,\dots,v_n,\beta)$ as in the statement of
Theorem~\ref{car-fibrant-lscts}.  Since $X$ is
$\mathcal{S}$-injective, the labelling map $\mu$ is one-to-one by
Proposition~\ref{simpl}. Therefore $u_i=v_i$ for $1\leq i \leq n$.
\epf

In particular, all cubical transition systems of the form
$\bl^{\cts}_{\mathcal{S}}(X)$ and all regular transition systems of
the form $\bl^{\rts}_{\mathcal{S}}(X)$ are fibrant because they are
$\mathcal{S}$-injective.

\appendix

\section{Proof of Proposition~\ref{compare-cellS}}
\label{cellS}

\bp \label{px-cts} Let $x\in \Sigma$. Every pushout of $p_x:C_1[x]
\sqcup C_1[x] \to \dd{x}$ in $\cts$ is bijective on states, and onto
on actions. There exists a pushout of $p_x$ which is not onto on
transitions. \ep

\bpf The category $\cts$ is a full coreflective category of $\wts$,
which means that the colimits in $\cts$ are calculated in
$\wts$. Therefore the forgetful functors taking a cubical transition
system to their sets of states and actions are
colimit-preserving. Since $p_x$ is bijective on states (onto on
actions resp.), any pushout of $p_x$ in $\cts$ is therefore bijective
on states (onto on actions resp.).

Let $x\in \Sigma$. One has $\omega (C_3[x,x,x]) =
(\{0,1\}^3,\{(x,1),(x,2),(x,3)\})$ by Proposition~\ref{cas_cube}.
Consider the quotient set
\[S = \{0,1\}^3\p\{-,+\}/\left( (0,0,0,-)=(0,0,0,+)=I\hbox{ and
  }(1,1,1,-)=(1,1,1,+)=F\right).\] Let 
\[W=(S,\{u_1,u^0,u^1,u_3\})\in \set^{\{s\}\cup \Sigma}\] with
$\mu(u_1)=\mu(u^0)=\mu(u^1)=\mu(u_3)=x$. For $\alpha\in \{-,+\}$,
consider the map \[\phi^\alpha:\omega (C_3[x,x,x])\to W\] of
$\set^{\{s\}\cup \Sigma}$ induced by the mappings
$(\epsilon_1,\epsilon_2,\epsilon_3)\mapsto
(\epsilon_1,\epsilon_2,\epsilon_3,\alpha)$ for
$(\epsilon_1,\epsilon_2,\epsilon_3)\in \{0,1\}^3$, $(x,1)\mapsto u_1$,
$(x,2)\mapsto u^\alpha $ and $(x,3)\mapsto u_3$.  Consider the
$\omega$-final lift $\overline{W}$ of the cone of maps
$\phi^-,\phi^+:\omega (C_3[x,x,x]) \rightrightarrows W$. By
Theorem~\ref{cubical-lift}, the weak transition system $\overline{W}$
is cubical.  Finally, consider the pushout diagram of cubical
transition systems:
\[
\xymatrix
{
C_1[x]\sqcup C_1[x] \fR{} \fD{p_{x}} && \overline{W} \fD{} \\
&& \\
\dd{x} \fR{} && \overline{\overline{W}} \cocartesien
}
\]
where the top horizontal arrow sends the $1$-transition $(0,(x,1),1)$
of the left-hand copy of $C_1[x]$ to $((1,0,0,-),u^-,(1,1,0,-))$ and
the $1$-transition $(0,(x,1),1)$ of the right-hand copy of $C_1[x]$ to
$((1,0,0,+),u^+,(1,1,0,+))$. We claim that the map of cubical
transition systems
\[\overline{W}\longrightarrow \overline{\overline{W}}\] is not surjective on
transitions. Indeed $\overline{W}$ contains the transitions
$(I,u_1,u^\alpha,u_3,F)$ for $\alpha\in \{-,+\}$, and the four transitions
\[(I,u_1,(1,0,0,-)),((1,0,0,-),u^-,u_3,F),(I,u_1,u^+,(1,1,0,+)),((1,1,0,+),u_3,F).\]
The cubical transition system $\overline{W}$ does not contain any
transition from $(1,0,0,-)$ to $(1,1,0,+)$. In the pushout
$\overline{\overline{W}}$, the identification $u^-=u^+$ is
made. Therefore from the five preceding transitions, one obtains by
using the composition axiom a transition $((1,0,0,-),u^-,(1,1,0,+))$.
\epf

\bp \label{surjectivity-cellS-cts} Every map of
$\cell_{\cts}(\mathcal{S})$ is bijective on states and onto on
actions. There exists a map of $\cell_{\cts}(\mathcal{S})$ which is
not onto on transitions.  \ep

\bpf A map of cubical transition systems is onto on actions if and
only if it satisfies the RLP with respect to the maps $\varnothing\to
\underline{x}$ for any $x\in \Sigma$. As a consequence, the class of
maps of cubical transition systems which are onto on actions is
accessible and accessibly embedded in the category of maps of cubical
transition systems by \cite[Proposition~3.3]{MR2506258}. Hence any map
of $\cell_{\cts}(\mathcal{S})$ is onto on actions. All maps of
$\mathcal{S}$ are bijective on states. Since the state set functor
from $\cts$ to $\set$ is colimit-preserving, all maps of
$\cell_{\cts}(\mathcal{S})$ are bijective on states.  The last
assertion is a corollary of Proposition~\ref{px-cts}.  \epf

\bp \label{px-onto} Let $x\in \Sigma$. Every pushout of $p_x:C_1[x]
\sqcup C_1[x] \to \dd{x}$ in $\rts$ is onto on states, on actions and
on transitions. \ep

\bpf Consider a pushout diagram in $\rts$: 
\[
\xymatrix
{
C_1[x] \sqcup C_1[x] \fD{p_x} \fR{} && X \fD{f} \\
&& \\
\dd{x} \fR{} && X'. \cocartesien 
}
\]
The category $\rts$ is a full reflective subcategory of
$\cts$. Therefore a colimit in $\rts$ is calculated by taking the
image by the reflection $\CSA_2:\cts\to \rts$ of the colimit in
$\cts$. The canonical map $Z\to \CSA_2(Z)$ is onto on states and
bijective on actions for all cubical transition systems $Z$ by
Proposition~\ref{calcul_CSA2}.  Therefore by Proposition~\ref{px-cts},
the map $f$ is onto both on states and on actions.  Let $X=(S,\mu:L\to
\Sigma,T)$ and $X'=(S',\mu:L'\to \Sigma,T')$. Write $f(T)$ for the set
of transitions of $X'$ of the form
$(f(\alpha),f(u_1),\dots,f(u_n),f(\beta))$ such that the tuple
$(\alpha,u_1,\dots,u_n,\beta)$ belongs to $T$. One has $f(T)\subset
T'$. Let $f(u)$ be an action of $X'$. Then there exists a transition
$(\alpha,u,\beta)$ of $X$ since $X$ is cubical. Therefore the tuple
$(f(\alpha),f(u),f(\beta))$ belongs to $f(T)$. This means that all
actions of $X'$ are used by a transition of $f(T)$.  Let
$(f(\alpha),f(u_1),\dots,f(u_n),f(\beta))$ be a transition of $f(T)$.
Then $(\alpha,u_{\sigma(1)},\dots,u_{\sigma(n)},\beta)$ is a
transition of $X$ for all permutations $\sigma$ of $\{1,\dots,n\}$. So
the tuple
$(f(\alpha),f(u_{\sigma(1)}),\dots,f(u_{\sigma(n)}),f(\beta))$ is a
transition of $f(T)$. Let $n\geq 3$ and $p,q \geq 1$ with $p+q<n$. Let
  \begin{multline*}(\alpha,u_1, \dots, u_n, \beta),(\alpha,u_1, \dots,
u_p, \mu),(\mu, u_{p+1}, \dots, u_n, \beta),\\(\alpha, u_1, \dots,
u_{p+q}, \nu),(\nu, u_{p+q+1}, \dots, u_n, \beta)
\end{multline*} be five transitions of $f(T)$. Let $(\alpha,u_1,
\dots, u_n, \beta) = (f(\gamma),f(v_1),\dots,f(v_n),f(\delta))$. There
exist two states $\epsilon$ and $\eta$ of $X$ such that the five
tuples $(\gamma,v_1,\dots,v_p,\epsilon)$,
$(\gamma,v_1,\dots,v_{p+q},\eta)$,
$(\epsilon,v_{p+1},\dots,v_n,\delta)$,
$(\eta,v_{p+q+1},\dots,v_n,\delta)$ and
$(\epsilon,v_{p+1},\dots,v_{p+q},\eta)$ are transitions of $X$ since
$X$ is cubical and by using the composition axiom in $X$. Therefore,
the five tuples
\begin{multline*}(f(\gamma),f(v_1),\dots,f(v_p),f(\epsilon)),
(f(\gamma),f(v_1),\dots,f(v_{p+q}),f(\eta)),\\
(f(\epsilon),f(v_{p+1}),\dots,f(v_n),f(\delta)),
(f(\eta),f(v_{p+q+1}),\dots,f(v_n),f(\delta)),\\
(f(\epsilon),f(v_{p+1}),\dots,f(v_{p+q}),f(\eta))
\end{multline*}
 are transitions of
$f(T)$. So the five tuples 
\begin{multline*}
(\alpha,u_1,\dots,u_p,f(\epsilon)),
(\alpha,u_1,\dots,u_{p+q},f(\eta)),\\
(f(\epsilon),u_{p+1},\dots,u_n,\beta),
(f(\eta),u_{p+q+1},\dots,u_n,\beta),\\
(f(\epsilon),u_{p+1},\dots,u_{p+q},f(\eta))
\end{multline*}
are transitions of $f(T)$. The point is that $X'$ is regular. One
deduces $f(\epsilon)=\mu$ and $f(\eta)=\nu$. One obtains
$(\mu,u_{p+1},\dots,u_{p+q},\nu)=(f(\epsilon),f(v_{p+1}),\dots,f(v_{p+q}),f(\eta))
\in f(T)$.  Let $n\geq 2$ and $1\leq p<n$. Let
$(f(\alpha),f(u_1),\dots,f(u_n),f(\beta))$ be a transition of
$f(T)$. Since $X$ is cubical, there exists a state $\mu$ such that
$(\alpha,u_1,\dots,u_p,\mu)$ and $(\mu,u_{p+1},\dots,u_n,\beta)$ are
two transitions of $X$. Since $X'$ is cubical, there exists a state
$\nu$ of $X'$ such that $(f(\alpha),f(u_1),\dots,f(u_p),\nu)$ and
$(\nu,f(u_{p+1}),\dots,f(u_n),f(\beta))$ are transitions of
$X'$. Since $X'$ is regular, one has $f(\mu)=\nu$.  Therefore
$(f(\alpha),f(u_1),\dots,f(u_p),\nu)$ and
$(\nu,f(u_{p+1}),\dots,f(u_n),f(\beta))$ belong to $f(T)$.  We have
proved that the tuple $Y=(S',L'\to \Sigma,f(T))$ is a regular
transition system. The map $X\to X'$ factors uniquely as a composite
$X\to Y \to X'$.  The map $\dd{x}\to X'$ factors uniquely as a
composite $\dd{x}\to Y \to X'$. By the universal property of the
pushout, one obtains $X'=Y$ and $T'=f(T)$.  \epf

\bp \label{states-acc} A map of regular transition systems is onto on
states if and only if it satisfies the RLP with respect to the map
$\varnothing\to \{0\}$. The class of maps of regular transition
systems which are onto on states is accessible and accessibly embedded
in the category of maps of regular transition systems. \ep

\bpf The first assertion is obvious. The second assertion is then a
consequence of \cite[Proposition~3.3]{MR2506258}.  \epf

\bp \label{trans-acc} A map of regular transition systems is onto on
transitions if and only if it satisfies the RLP with respect to the
maps $\varnothing\to C_n[x_1,\dots,x_n]$ for $n\geq 1$ and
$x_1,\dots,x_n\in \Sigma$. The class of maps of regular transition
systems which are onto on transitions is accessible and accessibly
embedded in the category of maps of regular transition systems. \ep

\bpf Let $f:X\to Y$ be a map of regular transition systems which is
onto on transitions. Consider a commutative diagram of weak transition
systems with $X$ and $Y$ regular:
\[
\xymatrix
{
&& C_n[x_1,\dots,x_n] \fR{k_2} && X \fD{f} \\
&& &&\\
\ar@{-->}[rruu]^-{k_1}\ar@{-->}[rrrruu]^-{\ell}C_n[x_1,\dots,x_n]^{ext} \fR{\subset}&& C_n[x_1,\dots,x_n] \fR{\phi} && Y.
}
\]
The lift $\ell$ exists since the map $f:X\to Y$ is onto on transitions
by hypothesis. Since $X$ is cubical, the map $\ell:
C_n[x_1,\dots,x_n]^{ext}\to X$ factors as a composite \[\ell:
C_n[x_1,\dots,x_n]^{ext}\stackrel{k_1}\longrightarrow
C_n[x_1,\dots,x_n] \stackrel{k_2}\longrightarrow X.\] The point is
that $Y$ is regular. Thus, $Y$ is orthogonal to the inclusion
$C_n[x_1,\dots,x_n]^{ext} \subset C_n[x_1,\dots,x_n]$ by
\cite[Theorem~5.6]{hdts}. Therefore $k_1$ is the inclusion
$C_n[x_1,\dots,x_n]^{ext} \subset C_n[x_1,\dots,x_n]$ and $\phi=f\circ
k_2$. We deduce that $f$ satisfies the RLP with respect to the maps
$\varnothing\to C_n[x_1,\dots,x_n]$ for $n\geq 1$ and
$x_1,\dots,x_n\in \Sigma$.

Conversely, let us suppose that $f:X\to Y$ is a map of regular
transition systems which satisfies the RLP with respect to the maps
$\varnothing\to C_n[x_1,\dots,x_n]$ for $n\geq 1$ and
$x_1,\dots,x_n\in \Sigma$. Let $(\alpha,u_1,\dots,u_n,\beta)$ be a
transition of $Y$. It yields a map
$C_n[\mu(u_1),\dots,\mu(u_n)]^{ext}\to Y$. Since $Y$ is cubical, this
map factors as a composite \[C_n[\mu(u_1),\dots,\mu(u_n)]^{ext}\subset
C_n[\mu(u_1),\dots,\mu(u_n)]\to Y.\] By hypothesis, the right-hand map
factors as a composite $C_n[\mu(u_1),\dots,\mu(u_n)]\to
X\stackrel{f}\to Y$.  Thus, the map
$C_n[\mu(u_1),\dots,\mu(u_n)]^{ext}\to Y$ factors as a composite
\[C_n[\mu(u_1),\dots,\mu(u_n)]^{ext}\to X \to Y.\] Hence $f$ is onto on
transitions.

The last assertion is then a consequence of
\cite[Proposition~3.3]{MR2506258}.  \epf

\bp \label{surjectivity-cellS} Every map of
$\cell_{\rts}(\mathcal{S})$ is onto on states, on actions and on
transitions.  \ep

\bpf A map of $\cell_{\rts}(\mathcal{S})$ is a transfinite composition
of maps which are onto on states and on transitions by
Proposition~\ref{px-onto}. By Proposition~\ref{states-acc} and
Proposition~\ref{trans-acc}, every map of $\cell_{\rts}(\mathcal{S})$
is then onto on states and on transitions. Let $f:X\to Y$ be a map of
$\cell_{\rts}(\mathcal{S})$. Let $u$ be an action of $Y$. Then there
exists a transition $(\alpha,u,\beta)$ of $Y$. Consider the map
$C_1[\mu(u)]\to Y$ taking the $1$-transition of $C_1[\mu(u)]$ to
$(\alpha,u,\beta)$. Then it factors as a composite $C_1[\mu(u)]\to X
\to Y$. The image of the $1$-transition of $C_1[\mu(u)]$ by the
left-hand map yields a $1$-transition $(\gamma,v,\delta)$ of $X$ such
that $(f(\gamma),f(v),f(\delta))=(\alpha,u,\beta)$. Therefore $f(v)=u$
and $f$ is onto on actions.  \epf

\section{Restricting an adjunction to a full reflective subcategory}
\label{restriction-adjunction-reflective-subcat}

The following proposition provides a tool to easily restrict the
cylinder and the path functors of cubical transition systems to the
reflective subcategory of regular ones. It is stated in a more general
setting than the one of locally presentable categories.

\bp \label{reflective-adjunction} Let $\LL\subset \K$ be two
categories with $\LL$ full and reflective. Let $R:\K\to \LL$ be the
reflection. Consider an adjunction $F \dashv G : \K\to \K$. Then the
following conditions are equivalent:
\begin{itemize}
\item[(i)] $F(\LL)\subset \LL$ and $G(\LL)\subset \LL$.
\item[(ii)] There is a natural isomorphism $R(F(X)) \iso F(R(X))$ for
  every $X \in \K$. 
\end{itemize}
If one of the two preceding conditions is satisfied, the restriction
of $F$ to $\LL$ is left adjoint to the restriction of $G$ to $\LL$.
\ep

\bpf The last assertion easily follows from the sequence of
isomorphisms
\[\LL(F(A),B)\iso \K(F(A),B)\iso \K(A,G(B))\iso \LL(A,G(B))\]
for any $A,B\in \LL$ and from the fact that $\LL$ is a full
subcategory of $\K$.

Let us prove now the implication $(i)\Rightarrow (ii)$. For any object
$X$ of $\K$ and any object $A$ of $\LL$, one has:
\begin{align*}
\LL(R(F(X)),A) &\iso \K(F(X),A) & \hbox{because $R$ is the left adjoint of $\LL\subset \K$}\\
&\iso \K(X,G(A)) & \hbox{because $G$ is the right adjoint of $F$}\\
&\iso \LL(R(X),G(A))& \hbox{by adjunction and since $G(A)\in \LL$}\\
&\iso \K(R(X),G(A))& \hbox{because $\LL$ is a full subcategory of $\K$}\\
&\iso \K(F(R(X)),A) & \hbox{because $G$ is the right adjoint of $F$}\\
&\iso \LL(F(R(X)),A)& \hbox{because $\LL$ is full in  $\K$ and $F(\LL)\subset \LL$.}
\end{align*}
By Yoneda applied in $\LL$, one obtains the natural isomorphism
$R(F(X)) \iso F(R(X))$.

Let us prove now the implication $(ii)\Rightarrow (i)$. Let $A$ be an
object of $\LL$. Then the unit map $\eta_A:A\to R(A)$, which is an
isomorphism since $A\in \LL$, gives rise to the isomorphism $F(A)\iso
F(R(A))$.  By $(ii)$, one then obtains the isomorphism $F(A)\iso
R(F(A))$. Hence $F(A)\in \LL$. We want to prove now that $G(A)\in
\LL$. One has the sequence of bijections
\begin{align*}
\K(G(A),G(A)) &\iso \K(F(G(A)),A)& \hbox{because $G$ is the right adjoint of $F$}\\
&\iso \LL(R(F(G(A)),A)& \hbox{because $R$ is the left adjoint of $\LL\subset \K$}\\
&\iso \K(R(F(G(A)),A)& \hbox{because $\LL$ is a full subcategory of $\K$}\\
&\iso \K(F(R(G(A)),A)& \hbox{because of $(ii)$}\\
&\iso \K(R(G(A)),G(A))& \hbox{because $G$ is the right adjoint of $F$.}
\end{align*} 
This means that the identity of $G(A)$ factors as a composite
\[G(A)\stackrel{\eta_{G(A)}}\longrightarrow R(G(A))
\stackrel{r}\longrightarrow G(A),\] i.e $r\circ \eta_{G(A)} =
\id_{G(A)}$. Hence $\eta_{G(A)}$ has a left inverse. We follow now the
argument of \cite{MO_reflective}.  By using the naturality of the unit
$\eta:\id\to R$, one obtains the commutative diagram
\[
\xymatrix{
R(G(A)) \fR{r} \fD{\eta_{R(G(A))}} && G(A) \fD{\eta_{G(A)}} \\
&& \\
R(R(G(A)) \fR{Rr} && R(G(A)).
}
\]
Since $r\circ \eta_{G(A)} =
\id_{G(A)}$, one has 
\[Rr \circ R(\eta_{G(A)}) = R(r\circ \eta_{G(A)}) = R(\id_{G(A)}) =
\id_{R(G(A))}.\] For all objects $Z$ of $\K$, the map
$R(\eta_Z):R(Z)\to R(R(Z))$ is an isomorphism by the universal
property of the reflection $R$.  With $Z= G(A)$, one obtains that
$R(\eta_{G(A)})$ is an isomorphism. Therefore $Rr =
R(\eta_{G(A)})^{-1}$ is an isomorphism. The map $\eta_{R(G(A))}$ is an
isomorphism as well since $\eta_{R(G(A))}=R(\eta_{G(A)})$. Therefore
\[\eta_{G(A)}\circ \left( r \circ (Rr \circ \eta_{R(G(A))})^{-1}\right) = \id_{R(G(A))}.\]
Hence $\eta_{G(A)}$ has a right inverse. Thus, $\eta_{G(A)}:G(A)\to
R(G(A))$ is an isomorphism. Hence $G(A)\in \LL$.  \epf


\begin{thebibliography}{Gau14b}

\bibitem[AHS06]{topologicalcat}
J.~Ad{\'a}mek, H.~Herrlich, and G.~E. Strecker.
\newblock Abstract and concrete categories: the joy of cats.
\newblock {\em Repr. Theory Appl. Categ.}, (17):1--507 (electronic), 2006.
\newblock Reprint of the 1990 original [Wiley, New York; MR1051419].

\bibitem[AR94]{MR95j:18001}
J.~Ad{\'a}mek and J.~Rosick{\'y}.
\newblock {\em Locally presentable and accessible categories}.
\newblock Cambridge University Press, Cambridge, 1994.

\bibitem[Bek00]{MR1780498}
T.~Beke.
\newblock Sheafifiable homotopy model categories.
\newblock {\em Math. Proc. Cambridge Philos. Soc.}, 129(3):447--475, 2000.

\bibitem[CS96]{MR1461821}
G.~L. Cattani and V.~Sassone.
\newblock Higher-dimensional transition systems.
\newblock In {\em 11th Annual IEEE Symposium on Logic in Computer Science (New
  Brunswick, NJ, 1996)}, pages 55--62. IEEE Comput. Soc. Press, Los Alamitos,
  CA, 1996.

\bibitem[Gau08]{ccsprecub}
P.~Gaucher.
\newblock Towards a homotopy theory of process algebra.
\newblock {\em Homology Homotopy Appl.}, 10(1):353--388 (electronic), 2008.

\bibitem[Gau10a]{symcub}
P.~Gaucher.
\newblock Combinatorics of labelling in higher dimensional automata.
\newblock {\em Theoretical Computer Science}, 411(11-13):1452--1483, 2010.
\newblock doi:10.1016/j.tcs.2009.11.013.

\bibitem[Gau10b]{hdts}
P.~Gaucher.
\newblock Directed algebraic topology and higher dimensional transition
  systems.
\newblock {\em New York J. Math.}, 16:409--461 (electronic), 2010.

\bibitem[Gau11]{cubicalhdts}
P.~Gaucher.
\newblock Towards a homotopy theory of higher dimensional transition systems.
\newblock {\em Theory Appl. Categ.}, 25:No.\ 25, 295--341 (electronic), 2011.

\bibitem[Gau14a]{erratum_cubicalhdts}
P.~Gaucher.
\newblock Erratum to "towards a homotopy theory of higher dimensional
  transition systems".
\newblock {\em Theory Appl. Categ.}, 29:No.\ 2, 17--20 (electronic), 2014.

\bibitem[Gau14b]{homotopyprecubical}
P.~Gaucher.
\newblock Homotopy theory of labelled symmetric precubical sets.
\newblock {\em New York J. Math.}, 20:93--131 (electronic), 2014.

\bibitem[Hir03]{ref_model2}
P.~S. Hirschhorn.
\newblock {\em Model categories and their localizations}, volume~99 of {\em
  Mathematical Surveys and Monographs}.
\newblock American Mathematical Society, Providence, RI, 2003.

\bibitem[Hov99]{MR99h:55031}
M.~Hovey.
\newblock {\em Model categories}.
\newblock American Mathematical Society, Providence, RI, 1999.

\bibitem[KR05]{ideeloc}
A.~Kurz and J.~Rosick{\'y}.
\newblock Weak factorizations, fractions and homotopies.
\newblock {\em Applied Categorical Structures}, 13(2):pp.141--160, 2005.

\bibitem[Low14]{MO_reflective}
Z.~L. Low.
\newblock About reflective full subcategories and small-orthogonality classes.
\newblock MathOverflow, 2014.
\newblock URL:http://mathoverflow.net/q/161463 (version: 2014-03-26).

\bibitem[Ols09a]{MO}
M.~Olschok.
\newblock Left determined model structures for locally presentable categories.
\newblock {\em Applied Categorical Structures}, 2009.
\newblock 10.1007/s10485-009-9207-2.

\bibitem[Ols09b]{MOPHD}
M.~Olschok.
\newblock {\em On constructions of left determined model structures}.
\newblock PhD thesis, Masaryk University, Faculty of Science, 2009.

\bibitem[Rap09]{1185.18014}
G.~Raptis.
\newblock {On the cofibrant generation of model categories.}
\newblock {\em J. Homotopy Relat. Struct.}, 4(1):245--253, 2009.

\bibitem[Ros09]{MR2506258}
J.~Rosick{\'y}.
\newblock On combinatorial model categories.
\newblock {\em Appl. Categ. Structures}, 17(3):303--316, 2009.

\bibitem[RR15]{rankweak}
G.~Raptis and J.~Rosick\'y.
\newblock The accessibility rank of weak equivalences.
\newblock {\em Theory Appl. Categ.}, 30:No.\ 19, 687--703 (electronic), 2015.

\end{thebibliography}

\end{document}